\documentclass[12pt]{extarticle}
\usepackage{amsfonts,amsmath,amsxtra,amsthm,amssymb,latexsym}
\usepackage[cp1251]{inputenc}
\usepackage{color}
\usepackage{indentfirst}
\usepackage{graphicx}

\usepackage{algorithm}
\usepackage{algorithmic}

\makeatletter
\@addtoreset{equation}{section} 
\makeatother

\numberwithin{equation}{section}
\newtheorem{theorem}{Theorem}
\newtheorem{lemma}{Lemma}
\newtheorem{corollary}{Corollary}

\newtheorem{remark}{Remark}

\newcommand{\Rr}{\mathbb{R}}
\newcommand{\Nn}{\mathbb{N}}

\usepackage{tikz}  
	\usetikzlibrary{arrows}
        \usetikzlibrary{shapes.geometric}
        \usetikzlibrary{shapes.symbols}
        \usetikzlibrary{positioning}
        \usetikzlibrary{calc}
        \usetikzlibrary{matrix} 
        \usetikzlibrary{fit}
        \usetikzlibrary{backgrounds}

\usepackage{geometry}
\begin{document}
\definecolor{Red}{rgb}{0, 0, 0}
\begin{center}
{\large\bf On adaptive discretization schemes for the regularization of ill-posed problems with semiiterative methods}

\vspace {0,2in}

{\sc W. Erb$^{(\dag)}$ \hspace{3cm} E.V. Semenova$^{(*)}$}

\vspace {0,2in}

\begin{minipage}{14cm} \it $(\dag)$ Institute of Mathematics, University of L\"ubeck\\
\phantom{$(\dag)$} Ratzeburger Allee 160, 23562 L\"ubeck, Germany
\end{minipage}

\vspace {0,1in}

\begin{minipage}{14cm} \it $(*)$ Institute of Mathematics, National Academy of Sciences of Ukraine\\
\phantom{$(*)$} Tereshchenkivska Str. 3, 01601 Kiev, Ukraine
\end{minipage}

%\vspace {0,2in}
%
%\emph{ e-mail:  lebedeva@ipnet.kiev.ua}
\end{center}

\vspace {0,1in}

\begin{abstract}
In this paper we investigate an adaptive discretization strategy for ill-posed linear problems combined with a regularization from a class of semiiterative methods. We show that 
such a discretization approach in combination with a stopping criterion as the discrepancy principle or the balancing principle 
yields an order optimal regularization scheme and allows to reduce the computational costs. 
\end{abstract}

\vspace{0.2in}


\section{Introduction}
%\textcolor[rgb]{1.00,0.00,0.00}{
Cost-efficient discretization methods for linear ill-posed equations are well-studied at the present day. Already in the early monographs of A.~N. Tikhonov, V.~K. Ivanov and M.~M. Lavrentiev
(see for example \cite{TA}) devoted to the regularization of unstable equations different discretizations were studied by applying finite-dimensional approximation techniques to the input data. 
But only in the beginning of the nineties R. Plato and G.~M. Vainikko \cite{PlVai} established estimates on the minimal rank of finite-dimensional operators that guarantee an efficient discretization 
and preserve the accuracy of the method at the same time. After this work the question on how to reduce
the volume of necessary discrete information was studied in several works. A first answer was given in \cite{P}. In this work, as a discretization domain the author proposed the
so-called hyperbolic cross for a particular class of ill-posed problems.
This cost-efficient discretization approach was then extended to larger classes of ill-posed problems in \cite{PerSol1}, \cite{Sol} and several other works. In particular, it is shown in \cite{SolLev} that if the solution of the ill-posed linear equation is smooth enough a large class of regularization methods is more economic when using the hyperbolic cross discretization instead of the classical 
discretization on a whole rectangular domain.

However, the results of these articles also imply that there are some cases when the discretization with the hyperbolic cross doesn't give an advantage in comparison 
with the classical discretization. In the case that the smoothness of the solution is low the volume of necessary discrete information
is the same for both mentioned approaches. Nevertheless, it is shown in \cite{MaaPerRamSol} that an additional adaptive strategy for the hyperbolic cross discretization
allows to reduce the volume of discrete information even for a low smoothness of the solution. In \cite{MaaPerRamSol} such a technique was applied for the Tikhonov-Philips regularization.
For a particular semiiterative method, namely for the $1$-method, a similar adaptive discretization approach was considered in \cite{SolVol} by S.~G. Solodky and E.~A. Volynets.

The work presented in this article is an extension of the results in \cite{MaaPerRamSol} and \cite{SolVol}. We show that with a similar adaptive discretization strategy as in \cite{SolVol}
it is possible to obtain cost-efficient and order optimal regularization schemes for a general class of semiiterative methods. Compared to \cite{SolVol}, slight changes in the adaptive algorithm allow us 
to control the influence of the discretization in the approximation error and to get explicit constants in the error estimates. Moreover, beside the discrepancy principle we consider in this work also 
the balancing principle as a stopping rule for the adaptive algorithm.

In the following two sections we introduce all preliminary information regarding semiiterative regularization methods and 
projection methods for the hyperbolic cross discretization
of linear ill-posed equations. The main results of this article can then be found in Section \ref{section-main}. 
Here, we present two adaptive regularization schemes (Algorithm \ref{algorithm1} and Algorithm \ref{algorithm2}) for 
the solution of ill-posed linear equations that combine regularization with semiiterative methods and a cost-efficient discretization strategy based on the hyperbolic cross. The algorithms are stopped either with the discrepancy principle
of Morozov or with the balancing principle. In Theorem \ref{teorem1} and Theorem \ref{teorem2} we will show that both algorithms are order optimal regularization methods
if the smoothness of the solution is contained in a given range. 
In the last section some numerical tests are provided that confirm the theoretical results on the order optimality of the algorithms. 

%Also this approach was used for discretization of $1$-method for $\mu\in (0,1)$. In the our paper we state the aim to extend the mentioned discretization technique for general class of semiitarative problems with reducing the standard estimation of volume for discrete information in comparison with standard.

%In this paper we consider the case when smoothness parameter belongs to wide interval from $0$ to $\infty$ (apoateriori case). Our results extends either the interval  for smoothness in comparison with previous from \cite{SolVol} and class of suitable semiiterative methods. Note that in \cite{SolVol} was considered only so called $1$-method.  For the stopping rule for proposed algorithm applyes the discrepancy principle.  %  From previous investigation it is easy to note that in some situation the direct discretization on some domain is not efficient in all situation. Some time it is occurs that
%The different aspects of   was considered in many famous publication of A.N. Tikh

\subsection{Preliminaries}

We shortly introduce the theoretical setup of this article. The space $X$ denotes a Hilbert space with inner product $(\cdot, \cdot)$ and norm $\|x\|=\sqrt{(x,x)}$. 
We consider operator equations of the first kind, i.e.
\begin{equation}\label{eq_main}
Ax=f,
\end{equation}
where $A: X \to X$ is a linear compact operator with  $ {\rm Range} (A) \ne
\overline {{\rm Range } (A)}$ and $f \in {\rm Range} (A)$.

By $x^\dag$ we denote the solution of \eqref{eq_main} with minimal
norm in $X$ that satisfies a H\"older-type source condition, namely we suppose that
\begin{equation}\label{source_condition}
x^\dag\in M_{\mu,\rho}(A) = \left\{x \in X:\, x=|A|^\mu v, \|v\|\leq \rho \right\} .
\end{equation}
In this article, $\rho > 0$ is supposed to be known whereas the unknown smoothness parameter $\mu$
is supposed to belong to an interval $(0,\mu_0]$ with some $0 <\mu_0< \infty$. Further, we set $|A| =
(A^*A)^{1/2}$ with $A^*$ denoting the adjoint operator of $A$. Also, we suppose that instead of the right-hand side $f$ in (\ref{eq_main})
we have given some perturbation $f_\delta$, $\|f-f_\delta\|\leq\delta$ with known noise level $\delta$.

Finally, we introduce the class of operators under consideration. We denote by
${\cal H}^r$ the class of compact linear operators $A,$ $\|A\|\leq 1,$ such
that for any $m\in \Nn$ the following conditions hold:
$$
\|(I-P_m)A\|\leq m^{-r}, \quad \|A(I-P_m)\|\leq m^{-r},
$$
where $P_m$ is the orthoprojector on the linear span of the first $m$ elements of some orthonormal basis $E=\{e_k\}_{k=1}^\infty$ in the space $X$.

\section{Properties of semiiterative methods}

In order to obtain cost-effective iterative schemes for the solution of \eqref{eq_main}, it is useful to consider sequences of orthogonal polynomials and use its 
three term recurrence relation to generate the iterates (see \cite{Hanke}). 
Let therefore $P_k(x)$, $k \in \Nn_0$, denote monic polynomials of degree $k$ that are orthogonal with respect to a weight function supported on the interval $[-1,1]$.
Then, the following recursion formula is valid (cf. \cite[I. Theorem 4.1]{Chihara}): 
\begin{equation} \label{Pol_cond}
P_{k+1}(x) = (x - \alpha_k) P_k(x) - \beta_k P_{k-1}(x), \quad P_0(x) = 1, \quad P_1(x) = x - \alpha_0.
\end{equation}
The recursion coefficients $\alpha_k \in \Rr$ and $\beta_k > 0$ are uniquely determined and $P_k(1) > 0$ holds for all $k \in \Nn_0$. 
Then (see \cite[Algorithm 1]{Erb14} and \cite[Theorem 2.1]{Hanke}) a semiiterative method can be defined as
\begin{equation}\label{SimiItMeth}
\begin{gathered}
x_{k} = x_{k-1} + ((1-\alpha_k) \omega_{k}-1) ( x_{k-1} - x_{k-2}) + 2 \omega_{k} \;  A^*( f - A x_{k-1}), \quad k \geq 1,\\
x_{-1} = 0, \quad x_0 = 2 \omega_0 \; A^* f,
\end{gathered}
\end{equation}
with
\[ \omega_{k} = \frac{1}{1-\alpha_k - \beta_k \omega_{k-1} } \quad \text{for}\quad k \geq 1 \quad \text{and} \quad \omega_0 = \frac{1}{1-\alpha_0}.\]
The iterative method (\ref{SimiItMeth}) yields an approximate solution of \eqref{eq_main} that can be written as
\begin{equation} \label{AppSol}
x_k=g_k(A^*A)A^*f,
\end{equation}
where $g_k(\lambda)$ is a polynomial of degree $k-1$. 
Combining (\ref{AppSol}) and (\ref{eq_main}), we can write the difference $x^\dag - x_k$ as
$$
x^\dag - x_k = r_k(A^*A)x^\dag,
$$
with the residual polynomial $r_k(\lambda)=1-\lambda g_k(\lambda)$ of degree $k$.
For the scheme (\ref{SimiItMeth}) the residual polynomials have the form
\begin{equation}\label{Resid}
r_k(\lambda)=\frac{P_{k}(1-2\lambda)}{P_k(1)},
\end{equation}
with the monic polynomials $P_k(\lambda)$ satisfying  the recursion formula (\ref{Pol_cond}).

Taking into account that we are only given a perturbed right-hand side $f_\delta$, the semiiterative method (\ref{SimiItMeth}) yields
\begin{equation}\label{AppSol_pert}
x_k^\delta=g_k(A^*A)A^*f_\delta
\end{equation}
as approximate solution of equation (\ref{eq_main}).

For an optimal speed of convergence to the solution $x^\dag\in M_{\mu,\rho}$, $0 < \mu \leq \mu_0 < \infty$, a sufficient condition for the residual polynomials $r_k$ in (\ref{Resid}) 
is (see \cite[Section 6.2.]{EnglHankeNeubauer})

\begin{equation}\label{cond1}
\begin{array}{lrll}
(i)  & |r_k(\lambda)|& \leq \; \kappa_0 & \mbox{for all} \quad \lambda\in [0,1], \; k\in \Nn, \\
(ii) & |\lambda^\frac{\mu}{2} r_k(\lambda)| & \leq \; \frac{\kappa_\mu}{(k+1)^{\mu}} & \mbox{for all} \quad  0 < \mu \leq \mu_0,\; \lambda\in[0,1], \; k\in \mathbb{N},
\end{array}
\end{equation}
with positive constants $\kappa_\mu > 0$, $0 < \mu \leq \mu_0$, and $\kappa_0 \geq 1$. The number $\mu_0$ is called the qualification of the semiiterative method (\ref{SimiItMeth}).

\vspace{0,2in}

{\bf Example 1:}  The $\nu$-methods of Brakhage \cite{Brakhage} are semiiterative methods based on the monic Jacobi polynomials 
$P_k^{(2\nu - \frac12, -\frac12)}$. For $\nu > 0$, its residual polynomials are given as
\[ r_k(\lambda) = \frac{P_k^{(2\nu - \frac12, -\frac12)}(1-2\lambda)}{P_k^{(2\nu - \frac12, -\frac12)}(1)}.\]
The qualification of this method is $\mu_0 = 2 \nu$. It is well-known that the Jacobi polynomials satisfy (see \cite[Section 4]{Hanke})
\[|r_k(\lambda)| \leq |r_k(0)| = 1, \quad \text{and} \quad |\lambda^{\nu} r_k(\lambda)| \leq |r_k(1)| = \binom{k+2\nu}{k}^{-1}.\]
Therefore $\kappa_0 = 1$, and, if we assume that $2 \nu$ is an integer, we get
\[ |\lambda^{\nu} r_k(\lambda)| \leq \frac{(2 \nu)! k!}{(k+2\nu)!} \leq \frac{(2\nu)!}{(k+1)^{2\nu}}.\]
Hence, in this case condition (ii) in \eqref{cond1} is satisfied with $\kappa_{2\nu} = (2 \nu)!$ and $\kappa_\mu \leq \max(\kappa_0,\kappa_{2\nu}) = (2 \nu)!$ for all $0 < \mu \leq 2 \nu$. 
The last statement follows from \cite[Theorem 4.2]{Hanke}. \\

{\bf Example 2:}  The $\nu$-method with parameter $\nu = \frac12$ is called Chebyshev method of Nemirovskii and Polyak \cite{NemirovskiiPolyak}. 
It is based on the Chebyshev polynomials $P_k^{(\frac12, -\frac12)}$ of fourth kind 
and its qualification is $\mu_0 = 1$. For the constants $\kappa_\mu$ in \eqref{cond1} we have $\kappa_\mu = 1$ for all $0 \leq \mu \leq 1$. \\

{\bf Example 3:}  The $\nu$-method with parameter $\nu = 1$ was studied in \cite{SolVol} in combination with cost-efficient adaptive discretization schemes. 
This scheme has qualification $\mu_0 = 2$ and the constants $\kappa_\mu$ in \eqref{cond1} can be chosen as $\kappa_0 = 1$ and $\kappa_\mu = 2$ for all $0 < \mu \leq 2$. \\

{\bf Example 4:}  Modified $\nu$-methods based on co-dilated Jacobi polynomials with an additional dilation parameter are investigated in \cite{Erb14}. These modified 
schemes have the same qualification $\mu_0 = 2 \nu$ as the $\nu$-methods. However, for these schemes the constant $\kappa_0$ is in general larger than $1$. 

\vspace{0,2in}

For our further analysis we need the Markov inequality in the form
\begin{equation}\label{markov}
|p'_k(\lambda)|\leq 2 \kappa k^2,
\end{equation}
where $p_k(\lambda)$, $\lambda \in[0,1]$, is a polynomial of degree $k$ and $\kappa =\max\limits_{\lambda \in[0,1]} |p_k(\lambda)|$.
Using the Markov inequality (\ref{markov}) and the conditions (i) and (ii) in \eqref{cond1} it is easy to obtain the following inequalities.

\begin{lemma} \label{lemma1}
If the polynomials $r_k(\lambda)$, $k \in \Nn$, satisfy (i), then
\begin{eqnarray} 
\sup_{\lambda \in [0,1]} |g_k(\lambda)| & \leq & 2 \kappa_0 k^2, \label{Ing1a}\\
\sup_{\lambda \in [0,1]} |\sqrt{\lambda} g_k(\lambda)| & \leq &  2 \kappa_0 k, \label{Ing1b}\\
|r_k(\lambda)-r_k(\tau)| & \leq &  2 \kappa_0 k^2 |\lambda-\tau| \quad \text{for all $\lambda, \tau \in [0,1]$.} \label{Ing1}
\end{eqnarray}
Moreover, if condition (ii) is satisfied with $\mu_0 \geq 2$, then
\begin{equation}\label{Ing2}
|\lambda r_k(\lambda)-\tau r_k(\tau)|\leq  2 \kappa_2 |\lambda-\tau| \quad \text{for all $\lambda, \tau \in [0,1]$.}
\end{equation}
\end{lemma}

{\bf Proof.}
Using the Markov inequality (\ref{markov}), condition (i) and the mean value theorem, we immediately get for $\lambda \in [0,1]$ the inequality
\begin{align*}
|g_k(\lambda)| = \frac{1-r_k(\lambda)}{\lambda} \leq \sup_{\lambda \in [0,1]} |r_k'(\lambda)| \leq 2 \kappa_0 k^2.
\end{align*}
Moreover, we get the inequality
\begin{align*}
\lambda g_k(\lambda)^2 = |\lambda g_k(\lambda)||g_k(\lambda)| = |1-r_k(\lambda)||g_k(\lambda)| \leq (1 + \kappa_0) 2 \kappa_0 k^2 \leq 4 \kappa_0^2 k^2.
\end{align*}
These two inequalities imply directly \eqref{Ing1a} and \eqref{Ing1b}. Using the same ingredients, we also get \eqref{Ing1}:  
\begin{align*}
\frac{r_k(\lambda) - r_k(\tau)}{\lambda-\tau} \leq \sup_{\lambda \in [0,1]} |r_k'(\lambda)| \leq 2 \kappa_0 k^2, \quad \lambda \neq \tau.
\end{align*}
Finally we show (\ref{Ing2}) in the case $\lambda\neq\tau$ (the case $\lambda=\tau$ is evident).
Using again the mean value theorem and the Markov inequality (\ref{markov}), we have
$$
\left|\frac{\lambda r_k(\lambda)-\tau r_k(\tau)}{\lambda-\tau}\right|\leq \max\limits_{\lambda\in[0,1]}|(\lambda r_k(\lambda))'|
\leq 2 (k+1)^2\max\limits_{\lambda\in[0,1]}|\lambda r_k(\lambda)|.
$$
Since the qualification $\mu_0$ of the residual polynomial is larger than $2$, we obtain due to property (ii) in \eqref{cond1} the inequality
$$
\left|\frac{\lambda r_k(\lambda)-\tau r_k(\tau)}{\lambda-\tau}\right|\leq 2 \kappa_2.
$$
\qed
\begin{lemma} \label{lemma2}
If the properties (i) and (ii) are satisfied with $\mu_0 \geq 2$, then the inequalities 
\begin{equation}\label{Ing3}
|\lambda(r_k(\lambda)-r_k(\tau))|\leq (\kappa_0 + 2 \kappa_2)|\lambda-\tau|,
\end{equation}
and
\begin{equation}\label{Ing4}
\sqrt{\lambda}|r_k(\lambda)-r_k(\tau)|\leq  2 \kappa_0 \sqrt{ \textstyle \frac12 + \frac{\kappa_2}{\kappa_0}}\, k |\lambda-\tau| 
\end{equation}
hold for $\lambda, \tau\in [0,1]$.
\end{lemma}
{\bf Proof.}
By (\ref{Ing2}) and (i), we immediately obtain
$$
|\lambda (r_k(\lambda)-r_k(\tau))|\leq |\lambda r_k(\lambda)-\tau r_k (\tau)|+|\lambda-\tau||r_k(\tau)| \leq
(\kappa_0 + 2\kappa_2) |\lambda-\tau|.
$$
The inequality (\ref{Ing4}) follows from (\ref{Ing3}) and (\ref{Ing1}) since
$$
\lambda|r_k(\lambda)-r_k(\tau)|^2=\lambda|r_k(\lambda)-r_k(\tau)||r_k(\lambda)-r_k(\tau)|\leq 2 \kappa_0 (\kappa_0 + 2\kappa_2) k^2|\lambda-\tau|^2.
$$ \qed

We remark that the qualitative statements of Lemma \ref{lemma1} and \ref{lemma2} are not new. In a simplified version, they can for instance be 
found in \cite{SolVol} for the residual polynomials of the $1$-method. However, for the error estimates in the following sections the explicit constants on the 
right hand side of the above inequalities, in particular of \eqref{Ing4}, play an important role. Therefore, we decided to include both lemmas with proof in this article. 

\section{Discretization schemes for linear equations}

To obtain finite dimensional approximations $A_\Omega$ of the operator $A$ we consider in this article projection methods and use the inner products
\begin{equation}\label{inn_pr}
(Ae_i, e_j) , \, (f_\delta, e_j), \quad (i,j)\in\Omega
\end{equation}
as discrete information about the linear equation (\ref{eq_main}), where  $\Omega \subset \{(i,j) \in \Nn^2\}$ denotes some subdomain of the coordinate plane.

In the following, we denote by $R_{\Omega}^{\mu_0}$ the class of iterative methods that solve (\ref{eq_main}) approximatively using the discretization domain $\Omega$ 
and a chosen semiiterative method with qualification $\mu_0 > 0$. We denote the corresponding iterates by
\begin{equation}\label{AppSol_pert_disc}
x_{\Omega,k}^\delta := g_k(A^*_{\Omega} A_{\Omega})A^*_{\Omega} f_\delta,\qquad k \in \Nn.
 \end{equation}
%as approximate solution of (\ref{eq_main}) obtained by some semiiterative method as (\ref{}) with fixed $\nu$ for discretized operator $A_\Omega$ and the right-hand side $P_{2^2n}f_\delta.$
We have the following general error bound for the iterates $x_{\Omega,k}^\delta$.

\begin{lemma} \label{lemma3}
Let the solution $x^\dag\in M_{\mu,\rho}(A)$ of \eqref{eq_main} satisfy a H\"older-type source condition with smoothness $0 < \mu \leq \mu_0$. Then, 
for the accuracy of an iterative method in the class $R_{\Omega}^{\mu_0}$ we obtain the following estimate:
$$
\|x^\dag-x_{\Omega,k}^\delta\|\leq \kappa_\mu \rho k^{-\mu}+ 2 \kappa_0 k \delta + 2 \kappa_0 k^2 \|x^\dag\| \Big(\|A^*A-A^*_\Omega A_\Omega\|+\|A^*_\Omega(A_\Omega-A)\|\Big).
$$
\end{lemma}

{\bf Proof.} We decompose the total error in the three terms
$$
\|x^\dag- x_{\Omega,k}^\delta\|\leq\|x^\dag-x_k\|+\|g_k(A_\Omega A_\Omega)A_\Omega^*(f-f_\delta)\|+\|x_k-g_k(A^*_\Omega A_\Omega)A_\Omega^* f\|
$$
and estimate each summand separately. For the estimate of the first summand we use the fact that $x^\dag\in M_{\mu,\rho}(A)$ 
and adopt a well-known result for regularization filters (see \cite[Lemma 3.3.6]{Rieder}) with qualification $\mu_0 \geq \mu$. This gives
$$
\|x^\dag-x_k\|=  \kappa_\mu (k+1)^{-\mu} \rho \leq \kappa_\mu \rho k^{-\mu}.
$$
The second estimate can be easily deduced from (\ref{Ing1b}):
$$
\|g_k(A_\Omega A_\Omega)A_\Omega^*(f-f_\delta)\|\leq \sup_{\lambda\in [0,1]} |\lambda^{1/2} g_k(\lambda) |\|f-f_\delta\|\leq 2 \kappa_0 k\delta.
$$
Finally, we estimate the third summand using (\ref{Ing1a}) and (\ref{Ing1}):
\begin{align*}
\|x_k-g_k(A^*_\Omega A_\Omega)A_\Omega f\| & \leq \|(g_k(A^*A)A^*A-g_k(A^*_\Omega A_\Omega)A_\Omega^*A_\Omega)x^\dag\|+\|g_k(A^*_\Omega A_\Omega)A_\Omega^*(A_\Omega-A)x^\dag\| \\
                                           & \leq \|(r_k(A^*_\Omega A_\Omega)-r_k(A^*A))x^\dag\|+\|g_k(A^*_\Omega A_\Omega)\|\|A^*_\Omega( A_\Omega-A)x^\dag\| \\
                                           & \leq 2 \kappa_0 k^2\|x^\dag\|\|A^*_\Omega A_\Omega-A^*A\|+2 \kappa_0 k^2\|x^\dag\|\|A^*_\Omega( A_\Omega-A)\|.
\end{align*}
Combining these three estimates, we get the statement of the lemma. \qed \\

We remark that the qualitative statement of Lemma \ref{lemma3} can be already found in \cite[Lemma 3.3]{SolVol} (therein proven for the $1$-method). Since the 
constants in the error estimate are important for our purposes, we decided to formulate also the above proof explicitly, although it is quite similar to 
the proof in \cite[Lemma 3.3]{SolVol}. 

\bigskip

The efficiency of the discretization depends heavily on the structure of $\Omega$. 
One standard choice for $\Omega$ (suggested in \cite{PlVai}) is to take the rectangular domain $\Omega=[1,...,M]\times[1,...,N]$. In this case the discretized operator $A_\Omega$ has the form
$$
A_\Omega=P_{M} A P_{N}.
$$
A second more efficient possibility is to take the hyperbolic cross 
\begin{equation}\label{HyperCross}
\Gamma_n := \bigcup_{k=1}^{2n} (2^{k-1}, 2^k]\times [1, 2^{2n-k}] \cup \{1\}\times[1,2^{2n}] \subset \Nn^2
\end{equation}
as a discretization domain. In this case, the discretized operator $A_n := A_{\Gamma_n}$ has the form
\begin{equation}\label{A_disc}
A_n:= A_{\Gamma_n} = \sum_{k=1}^{2n}(P_{2^k}-P_{2^{k-1}})AP_{2^{2n-k}}+P_1AP_{2^{2n}}.
\end{equation}

\begin{figure}[h]
	\begin{center}
	\begin{tikzpicture}
		[scale=1]
	        \tikzset{myshade/.style={minimum size=0.4cm,shading=radial,inner color=white,outer color={#1!90!gray}}}
	        \tikzset{axe/.style={thick, ->, >=stealth'}}
		\tikzset{hc/.style={circle, thick, anchor = center, myshade = blue}}
		\tikzset{ze/.style={circle, thick, anchor = center, myshade = lightgray}}

		\matrix [matrix of math nodes, column sep={0.5cm,between origins}, row sep={0.5cm,between origins},
					] (A)
		{
              |[hc]| & |[hc]| & |[hc]| & |[hc]| & |[hc]| & |[hc]| & |[hc]| & |[hc]| & |[hc]| & |[hc]| & |[hc]| & |[hc]| & |[hc]| & |[hc]| & |[hc]| & |[hc]| & |[ze]| & |[ze]| & |[ze]| & |[ze]| \\
              |[hc]| & |[hc]| & |[hc]| & |[hc]| & |[hc]| & |[hc]| & |[hc]| & |[hc]| & |[ze]| & |[ze]| & |[ze]| & |[ze]| & |[ze]| & |[ze]| & |[ze]| & |[ze]| & |[ze]| & |[ze]| & |[ze]| & |[ze]| \\
              |[hc]| & |[hc]| & |[hc]| & |[hc]| & |[ze]| & |[ze]| & |[ze]| & |[ze]| & |[ze]| & |[ze]| & |[ze]| & |[ze]| & |[ze]| & |[ze]| & |[ze]| & |[ze]| & |[ze]| & |[ze]| & |[ze]| & |[ze]| \\
              |[hc]| & |[hc]| & |[hc]| & |[hc]| & |[ze]| & |[ze]| & |[ze]| & |[ze]| & |[ze]| & |[ze]| & |[ze]| & |[ze]| & |[ze]| & |[ze]| & |[ze]| & |[ze]| & |[ze]| & |[ze]| & |[ze]| & |[ze]| \\
              |[hc]| & |[hc]| & |[ze]| & |[ze]| & |[ze]| & |[ze]| & |[ze]| & |[ze]| & |[ze]| & |[ze]| & |[ze]| & |[ze]| & |[ze]| & |[ze]| & |[ze]| & |[ze]| & |[ze]| & |[ze]| & |[ze]| & |[ze]| \\
              |[hc]| & |[hc]| & |[ze]| & |[ze]| & |[ze]| & |[ze]| & |[ze]| & |[ze]| & |[ze]| & |[ze]| & |[ze]| & |[ze]| & |[ze]| & |[ze]| & |[ze]| & |[ze]| & |[ze]| & |[ze]| & |[ze]| & |[ze]| \\
              |[hc]| & |[hc]| & |[ze]| & |[ze]| & |[ze]| & |[ze]| & |[ze]| & |[ze]| & |[ze]| & |[ze]| & |[ze]| & |[ze]| & |[ze]| & |[ze]| & |[ze]| & |[ze]| & |[ze]| & |[ze]| & |[ze]| & |[ze]| \\
              |[hc]| & |[hc]| & |[ze]| & |[ze]| & |[ze]| & |[ze]| & |[ze]| & |[ze]| & |[ze]| & |[ze]| & |[ze]| & |[ze]| & |[ze]| & |[ze]| & |[ze]| & |[ze]| & |[ze]| & |[ze]| & |[ze]| & |[ze]| \\
              |[hc]| & |[ze]| & |[ze]| & |[ze]| & |[ze]| & |[ze]| & |[ze]| & |[ze]| & |[ze]| & |[ze]| & |[ze]| & |[ze]| & |[ze]| & |[ze]| & |[ze]| & |[ze]| & |[ze]| & |[ze]| & |[ze]| & |[ze]| \\
              |[hc]| & |[ze]| & |[ze]| & |[ze]| & |[ze]| & |[ze]| & |[ze]| & |[ze]| & |[ze]| & |[ze]| & |[ze]| & |[ze]| & |[ze]| & |[ze]| & |[ze]| & |[ze]| & |[ze]| & |[ze]| & |[ze]| & |[ze]| \\
              |[hc]| & |[ze]| & |[ze]| & |[ze]| & |[ze]| & |[ze]| & |[ze]| & |[ze]| & |[ze]| & |[ze]| & |[ze]| & |[ze]| & |[ze]| & |[ze]| & |[ze]| & |[ze]| & |[ze]| & |[ze]| & |[ze]| & |[ze]| \\
              |[hc]| & |[ze]| & |[ze]| & |[ze]| & |[ze]| & |[ze]| & |[ze]| & |[ze]| & |[ze]| & |[ze]| & |[ze]| & |[ze]| & |[ze]| & |[ze]| & |[ze]| & |[ze]| & |[ze]| & |[ze]| & |[ze]| & |[ze]| \\
              |[hc]| & |[ze]| & |[ze]| & |[ze]| & |[ze]| & |[ze]| & |[ze]| & |[ze]| & |[ze]| & |[ze]| & |[ze]| & |[ze]| & |[ze]| & |[ze]| & |[ze]| & |[ze]| & |[ze]| & |[ze]| & |[ze]| & |[ze]| \\
              |[hc]| & |[ze]| & |[ze]| & |[ze]| & |[ze]| & |[ze]| & |[ze]| & |[ze]| & |[ze]| & |[ze]| & |[ze]| & |[ze]| & |[ze]| & |[ze]| & |[ze]| & |[ze]| & |[ze]| & |[ze]| & |[ze]| & |[ze]| \\
              |[hc]| & |[ze]| & |[ze]| & |[ze]| & |[ze]| & |[ze]| & |[ze]| & |[ze]| & |[ze]| & |[ze]| & |[ze]| & |[ze]| & |[ze]| & |[ze]| & |[ze]| & |[ze]| & |[ze]| & |[ze]| & |[ze]| & |[ze]| \\
              |[hc]| & |[ze]| & |[ze]| & |[ze]| & |[ze]| & |[ze]| & |[ze]| & |[ze]| & |[ze]| & |[ze]| & |[ze]| & |[ze]| & |[ze]| & |[ze]| & |[ze]| & |[ze]| & |[ze]| & |[ze]| & |[ze]| & |[ze]| \\
              |[ze]| & |[ze]| & |[ze]| & |[ze]| & |[ze]| & |[ze]| & |[ze]| & |[ze]| & |[ze]| & |[ze]| & |[ze]| & |[ze]| & |[ze]| & |[ze]| & |[ze]| & |[ze]| & |[ze]| & |[ze]| & |[ze]| & |[ze]| \\
              |[ze]| & |[ze]| & |[ze]| & |[ze]| & |[ze]| & |[ze]| & |[ze]| & |[ze]| & |[ze]| & |[ze]| & |[ze]| & |[ze]| & |[ze]| & |[ze]| & |[ze]| & |[ze]| & |[ze]| & |[ze]| & |[ze]| & |[ze]| \\
              |[ze]| & |[ze]| & |[ze]| & |[ze]| & |[ze]| & |[ze]| & |[ze]| & |[ze]| & |[ze]| & |[ze]| & |[ze]| & |[ze]| & |[ze]| & |[ze]| & |[ze]| & |[ze]| & |[ze]| & |[ze]| & |[ze]| & |[ze]| \\
              |[ze]| & |[ze]| & |[ze]| & |[ze]| & |[ze]| & |[ze]| & |[ze]| & |[ze]| & |[ze]| & |[ze]| & |[ze]| & |[ze]| & |[ze]| & |[ze]| & |[ze]| & |[ze]| & |[ze]| & |[ze]| & |[ze]| & |[ze]| \\                        
                };
		% i - axis
		
		\node (x1) [above=10pt, label=above:{1}] at (A-1-1) {};
		\node (zero) [left = 10pt] at (x1) {};
		\node (x2) [above=10pt, label=above:{10}] at (A-1-10) {};
		\node (x3) [above=10pt, label=above:{20}] at (A-1-20) {};
		\node (x4) [right = 20 pt, label=above:{i}] at (x3) {};
		\draw[axe] (zero) +(-20pt,0) -- (x4.east) ;
		\draw (x1.north) -- (x1.south);
		\draw (x2.north) -- (x2.south);
		\draw (x3.north) -- (x3.south);
		% j - axis
		%\node (a) [right=40pt] at (A-1-7) {};
		%\node (b) [right=40pt] at (A-13-7) {};
		\node (y1) [left=10pt, label=left:{1}] at (A-1-1) {};
		\node (y2) [left=10pt, label=left:{10}] at (A-10-1) {};
		\node (y3) [left=10pt, label=left:{20}] at (A-20-1) {};
		\node (y4) [below = 20 pt, label=left:{j}] at (y3) {};
		\draw[axe] (zero) +(0,20pt) -- (y4.south) ;
		\draw (y1.west) -- (y1.east);
		\draw (y2.west) -- (y2.east);
		\draw (y3.west) -- (y3.east);
		
		\node [left=10pt, label=above:{$\Gamma_2$}] at (zero) {};
							
	\end{tikzpicture}
	\caption[Visualization of hyperbolic cross $\Gamma_2$]{Visualization of the hyperbolic cross $\Gamma_2$. The elements contained in $\Gamma_2$ are colored in blue. $\Gamma_2$ contains $48$ elements. 
	In contrast, the rectangular domain  $\Omega = [1, \ldots, 2^4] \times [1, \ldots, 2^4]$ contains $256$ elements.}
	\label{fig:arrangement-coefficients}
	\end{center}
\end{figure}
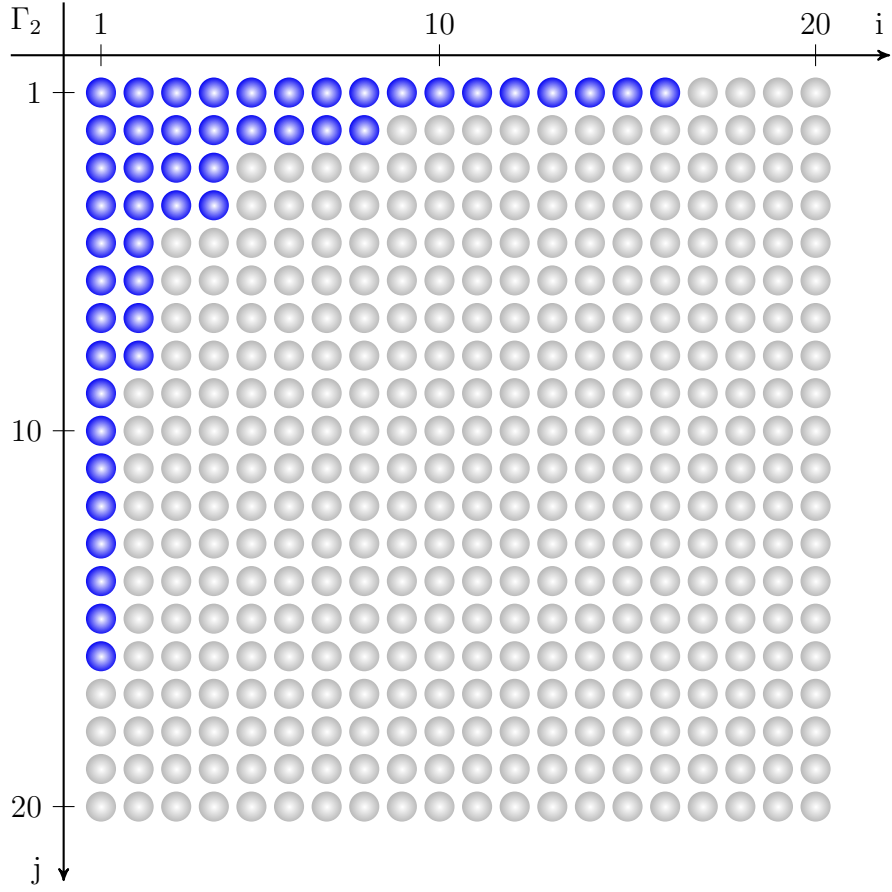

\begin{remark} It is easy to compute the volume of Galerkin information necessary for the realization of the hyperbolic cross approximation $A_n$ given by (\ref{A_disc}): one 
has to compute $\# \Gamma_n = 2^{2n}(n+1)$ inner products to construct $A_n$. In comparison, one has to compute $2^{4n}$ inner products for the approximation with $A_\Omega$ for
the standard quadratic domain $\Omega = [1, \ldots, N]^2$, $N = 2^{2n}$.
\end{remark}

If $A \in {\cal H}^r$ and $A_n$ has the form (\ref{A_disc}), the following error estimates hold (see \cite{Sol}):
\begin{equation}
\|A^*A-A^*_nA_n\|\leq (1+2^{r+3})2^{-2rn} n, \quad \|A_n^*(A-A_n)\|\leq 3 \; 2^{-2rn+r}n .
\end{equation}
Moreover, since $A \in {\cal H}^r$ we have the bound
\begin{equation}
\|A -A P_{2^{2n}}\| \leq 2^{-2rn}.
\end{equation}

\begin{corollary} \label{corollary-discretization}
If $A \in {\cal H}^r$, $A_n$ has the form (\ref{A_disc}) and 
\begin{equation} \label{conddiscr}
(1+2^{r+3})2^{-2rn} n < \frac{\gamma \delta}{2 k \rho }
\end{equation}
with $k \in \Nn$ and a control parameter $\gamma > 0$, then
\begin{equation} \label{inequality-discretization}
\|A^*A-A^*_nA_n\|  \leq \frac{\gamma \delta}{2 k \rho}, \qquad
\|A_n^*(A-A_n)\|  \leq \frac{\gamma \delta}{2 k \rho}, \qquad
\|A - A P_{2^{2n}} \|  \leq \frac{\gamma \delta}{2 k \rho n }.
\end{equation}
If $x^\dag \in M_{\mu,\rho}(A)$, the approximation error of methods from the class $R_{n}^{\mu_0} := R_{\Gamma_n}^{\mu_0}$ with $\mu_0 \geq \mu > 0$ is given by
\begin{equation}\label{accuracy}
\|x^\dag- x_{n,k}^\delta\|\leq \kappa_\mu \rho k^{-\mu}+2 \kappa_0 k \delta (1 + \gamma).
\end{equation}
\end{corollary}

In the adaptive algorithms of the next section the values $\delta$, $\rho$ and $k \in \Nn$ can not be chosen freely. The control parameter $\gamma$ in \eqref{conddiscr} enables 
a possible user to handle the trade-off between accuracy and cost-efficiency of the adaptive scheme.

\section{Regularization with semiiterative methods and adaptive discretization strategies} \label{section-main}

This section contains the main new results of the article. We present two algorithms in which the regularization with semiiterative methods is combined with
an adaptive and cost-efficient discretization strategy. As stopping rule for the regularization we consider the discrepancy principle of Morozov \cite{Morozov} on the one hand and
the balancing principle \cite{PerSch} on the other. 

\subsection{The discrepancy principle as stopping rule}

To solve the linear ill-posed problem (\ref{eq_main}) we consider the adaptive Algorithm \ref{algorithm1} which combines a method from the class $R_{n}^{\mu_0}$, $\mu_0 \geq 2$, with 
an adaptive discretization strategy and the discrepancy principle of Morozov \cite{Morozov} (implemented as \eqref{Discr_princp} in Algorithm \ref{algorithm1}) as stopping rule.

\begin{algorithm}
\caption{Adaptive algorithm to solve (\ref{eq_main}) using the discrepancy principle}
\label{algorithm1}

\begin{algorithmic}[H]
\vspace{2mm}
\STATE given data $A\in {\cal H}^r, \delta, f_\delta, \rho$.
\STATE choose control parameters $\gamma > 0$, $\tau > \kappa_0 \left( 1 + \sqrt{ \textstyle \frac12 + \frac{\kappa_2}{\kappa_0} } \, \gamma \right)$.
\STATE choose discretization level $n \geq 1$ such that $(1+2^{r+3})2^{-2rn} n < \frac{\gamma \delta}{2 \rho}$ holds. 
\STATE compute Galerkin information:
        $$
        (f_\delta, e_j), \quad j\in [1, 2^{2n}], \qquad
        (Ae_i, e_j), \quad (i,j) \in \Gamma_{n}.
        $$
\WHILE {(discrepancy principle $=$ false)}  \vspace{2mm}

\STATE choose $K_n \in \Nn$ as maximal integer such that 
        \begin{equation} \label{equation-n}
                (1+2^{r+3})2^{-2rn} n < \frac{\gamma \delta}{2 K_{n} \rho}
        \end{equation} 
        is satisfied. \\[2mm]
\FOR {($k = 1:K_n$)} \vspace{2mm}
\STATE compute $k-th.$ iterate of semiiterative method in the class $R_n^{\mu_0}$, $\mu_0 \geq 2$ (cf. \eqref{SimiItMeth}):
\STATE \[ x_{n,k}^\delta = x_{n,k-1}^\delta + ((1-\alpha_k) \omega_{k}-1) ( x_{n,k-1}^\delta - x_{n,k-2}^\delta) + 2 \omega_{k}  A_n^*( f_\delta - A_n x_{n,k-1}^\delta). \]
\IF {  \begin{equation} \label{Discr_princp}
        \|A_{n} x_{n,k}^\delta - P_{2^{2n}} f_\delta\|\leq \tau \delta \quad \text{and} \quad\|A_{n} x_{n,j}^\delta- P_{2^{2n}} f_\delta \| > \tau \delta \quad \text{for all $j < k$},
       \end{equation} } \vspace{2mm}
\STATE discrepancy principle $=$ true, \\[2mm]
\RETURN stopping index $K = k$, discretization level $n$ and solution $x_{n,K}^\delta$.\\[2mm]
\ENDIF \\[2mm]
\ENDFOR\\[2mm]
\STATE increase discretization level $n \to n+1$. 
\STATE compute new Galerkin information:
        $$
        (f_\delta, e_j), \quad j\in [2^{2n-2}, 2^{2n}]; \qquad
        (Ae_i, e_j), \quad (i,j) \in \Gamma_{n} \setminus \Gamma_{n-1}.
        $$
\ENDWHILE
\end{algorithmic}
\end{algorithm}

In Theorem \ref{teorem1} we show that Algorithm \ref{algorithm1} yields an order optimal regularization scheme for the solution of \eqref{eq_main}. 
For the proof we need some additional statements.

\begin{lemma}\label{lemma4}
Let $A \in {\cal H}^r$, $x^\dag \in M_{\mu,\rho}(A)$ and $x_{n,k}^\delta$ be computed according to Algorithm \ref{algorithm1}. Then, we have the estimate
$$
\|A x_k - f\|\leq \|A_n x_{n,k}^\delta - P_{2^{2n}} f_\delta\|+c_1\delta,
$$
with $c_1 = \kappa_0 + 2 + \left( \sqrt{ \kappa_0( \textstyle \frac{\kappa_0}2 + \kappa_2 )} + \frac{1}{2n} \right) \gamma.$
\end{lemma}

{\bf Proof.}
To prove the statement we write $Ax_k-f$ in the following telescoping sum:
$$
Ax_k-f = A x_k - A x_k^\delta + A x_k^\delta - A_n x_{n,k}^\delta + A_n x_{n,k}^\delta - P_{2^{2n}} f_\delta + P_{2^{2n}} f_\delta - P_{2^{2n}} f + P_{2^{2n}} f - f
$$
Now, using the triangle inequality and the definition of the generating polynomials $g_n$ we get
\begin{align*}
\| Ax_k-f \| & \leq  \|A g_k(A^*A) A^* (f-f_\delta)\| + \| A g_k(A^*A) A^* f_\delta - A_n g_k(A_n^* A_n) A_n^* f_\delta\| \\ 
& \quad + \|A_n x_{n,k}^\delta - P_{2^{2n}} f_\delta\| + \|P_{2^{2n}} (f_\delta - f)\| + \|P_{2^{2n}} f - f\|.
\end{align*}

Since $g_k(A^*A)A^*=A^*g_k(AA^*)$, we can further estimate
\begin{align*}
\| Ax_k-f \| & \leq  \sup_{\lambda\in[0,1]}|\lambda g_k(\lambda)| \|f-f_\delta\| + \| (r_k(AA^*) - r_k(A_n A_n^*)) Ax^\dag\| \\ 
& \quad + \|A_n x_{n,k}^\delta - P_{2^{2n}} f_\delta\| + \|f_\delta - f\| + \|(P_{2^{2n}}A - A) x^\dag\|.
\end{align*}

Using property (i) of the residual polynomials, inequality (\ref{Ing4}) as well as the inequalities \eqref{inequality-discretization} in Corollary \ref{corollary-discretization}, we finally obtain
\begin{align*}
\| Ax_k-f \| & \leq  ( \kappa_0 + 2 )\delta + 2 \sqrt{ \kappa_0( \textstyle \frac{\kappa_0}2 + \kappa_2 )}\, k \| A^*A - A_n^* A_n \| \rho + \|A_n x_{n,k}^\delta - P_{2^{2n}} f_\delta\| + \| P_{2^{2n}}A - A\| \rho. \\
             & \leq \|A_n x_{n,k}^\delta - P_{2^{2n}} f_\delta\| + \left( \kappa_0 + 2 + \left( \sqrt{ \kappa_0( \textstyle \frac{\kappa_0}2 + \kappa_2 )} + \frac{1}{2n} \right) \gamma \right) \delta.
\end{align*}
\qed

\begin{lemma}\label{lemma5}
Let $A \in {\cal H}^r$, $x^\dag \in M_{\mu,\rho}(A)$, $0 < \mu \leq \mu_0 - 1$, $\mu_0 \geq 2$ and $x_{n,k}^\delta$ be computed according to Algorithm \ref{algorithm1}. Further, we assume
that the level of noise satisfies $\delta < \|f\|$. Then, the stopping index $K$ is bounded by
\begin{equation} \label{inequalitystoppingindex}
K < c_2 \rho^{\frac{1}{\mu+1}} \delta^{-\frac{1}{\mu+1}}, \quad \text{with} 
\quad c_2 = \max \left\{ \left( \frac{\kappa_\mu}{\tau - \kappa_0 \left( 1 + \sqrt{ \textstyle \frac12 + \frac{\kappa_2}{\kappa_0} }\, \gamma\right) } \right)^{\frac{1}{\mu+1}} , 1 \right\}.
\end{equation}
If Algorithm \ref{algorithm1} is not stopped in the first iteration of the while loop, then also the index $K_{n-1}$ satisfies the estimate
\begin{equation} \label{inequalitymaximalindex}
K_{n-1}+1 < c_2 \rho^{\frac{1}{\mu+1}} \delta^{-\frac{1}{\mu+1}}.
\end{equation}
\end{lemma}

{\bf Proof.}
To prove this statement we consider the second inequality \eqref{Discr_princp} in Algorithm \ref{algorithm1}. For the stopping index $K \geq 2$, we get
\begin{align*}
\tau \delta & < \|A_{n} x_{n,K-1}^\delta- P_{2^{2n}} f_\delta \| = \|r_{K-1} (A_n A_n^*) f_\delta\| \\
              & \leq \|r_{K-1} (A A^*) A x^\dag\| + \|( r_{K-1} (A_n A_n^*) - r_{K-1} (A A^*)) A x^\dag \| + \|r_{K-1} (A_n A_n^*) (f-f_\delta)\|.
\end{align*}
Now using the fact that $x^\dag \in M_{\mu,\rho}(A)$, property (i) of the residual polynomials $r_n$ as well as inequality \eqref{Ing4} (here $\mu_0 \geq 2$ is necessary), we obtain
\[\tau \delta  < \||A|^{\mu+1} r_{K-1} (A^* A) v\| + 2 \sqrt{ \kappa_0( \textstyle \frac{\kappa_0}2 + \kappa_2 )}\, (K-1) \| A^*A - A_n^* A_n \| \rho + \kappa_0 \delta.\]
Finally, using property (ii) of \eqref{cond1} (here $\mu+1 \leq \mu_0$ must be satisfied) and Corollary \ref{corollary-discretization} (the conditions of the corollary
are satisfied by the construction of Algorithm \ref{algorithm1}) we conclude
\[\tau \delta  < \kappa_{\mu+1} K^{-(\mu+1)} \rho + \left(\kappa_0 + \sqrt{ \kappa_0( \textstyle \frac{\kappa_0}2 + \kappa_2 )}\, \gamma \right) \delta. \]
Solving this inequality for the index $K$ implies inequality \eqref{inequalitystoppingindex}. 

For the case that Algorithm \ref{algorithm1} is stopped at $K = 1$, the assumption $\delta < \|f\|$ implies
\[ \delta < \|f\| = \||A|^{\mu+1} v\| \leq \rho. \]
Thus, we get for $K = 1$
\[ K <  \rho^{\frac{1}{\mu+1}} \delta^{-\frac{1}{\mu+1}}.\]
The proof of inequality \eqref{inequalitymaximalindex} follows the same lines of argumentation as the proof of inequality \eqref{inequalitystoppingindex}
with $K-1$ replaced by $K_{n-1}$ and $n$ replaced by $n-1$. 
\qed

\begin{theorem}\label{teorem1}
Let $A \in {\cal H}^r$, $\delta < \|f\|$ and $\mu_0 \geq 2$ for the qualification of the semiiterative method. Then, Algorithm \ref{algorithm1} gives
an order optimal regularization method for the solution $x^\dag$ of \eqref{eq_main} in the class $M_{\mu,\rho}(A)$ for all $0 < \mu \leq \mu_0 - 1$. In particular,
the approximative solution $x_{n,K}^\delta$ given by Algorithm \ref{algorithm1} satisfies 
\begin{equation} \label{orderoptimalitydiscrepancy}
\|x^\dag-x_{n,K}^\delta\|\leq C \rho^{\frac{1}{\mu+1}} \delta^{\frac{\mu}{\mu+1}},
\end{equation}
with $C = \left(\kappa_0^{\frac{1}{\mu+1}} (\tau + c_1)^{\frac{\mu}{\mu+1}} + 2 \kappa_0 (1+\gamma) c_2 \right)$.
\end{theorem}

{\bf Proof.}
We only have to prove inequality \eqref{orderoptimalitydiscrepancy}, i.e. 
that Algorithm \ref{algorithm1} gives an order optimal reconstruction scheme for all $0 < \mu \leq \mu_0 -1$. Then, it follows from a general result of R. Plato \cite[Theorem 2.1]{Plato} that
Algorithm \ref{algorithm1} is also a regularization method for the solution of \eqref{eq_main}. 
Using the estimates of Lemma \ref{lemma3} and Corollary \ref{corollary-discretization} (excluding the estimate for the approximation error), we get the error bound.
\begin{equation} \label{ax1}
\|x^\dag- x_{n,K}^\delta\|\leq \| x^\dag - x_K \| + 2 \kappa_0 K \delta (1 + \gamma).
\end{equation}
To estimate the approximation error, we use an interpolation inequality (see \cite[Satz 2.4.2.]{Rieder} or \cite[Satz 2.3.3.]{Louis}) and obtain
\begin{align*}
\| x^\dag - x_K \| & \leq \|r_K(A^* A) x^\dag\| = \||A|^{\mu} r_K(A^* A) v\| \\
                   & \leq \||A|^{\mu+1} r_K(A^* A) v\|^{\frac{\mu}{\mu+1}} \| r_K(A^* A) v\|^{\frac{1}{\mu+1}} \\
                   & \leq \|A r_K(A^* A) x^\dag\|^{\frac{\mu}{\mu+1}} (\kappa_0 \rho)^{\frac{1}{\mu+1}} = \|A x_K - f\|^{\frac{\mu}{\mu+1}} (\kappa_0 \rho)^{\frac{1}{\mu+1}}.
\end{align*}
Now, Lemma \ref{lemma4} and the discrepancy principle \eqref{Discr_princp} give
\begin{align*}
\| x^\dag - x_K \| &\leq (\|A_n x_{n,k}^\delta - P_{2^{2n}} f_\delta\| + c_1\delta)^{\frac{\mu}{\mu+1}} (\kappa_0 \rho)^{\frac{1}{\mu+1}} \\
                   &\leq \kappa_0^{\frac{1}{\mu+1}} (\tau + c_1)^{\frac{\mu}{\mu+1}} \rho^{\frac{1}{\mu+1}}\delta^{\frac{\mu}{\mu+1}}.
\end{align*}
For the data error in \eqref{ax1} Lemma \ref{lemma5} gives the bound
\[2 \kappa_0 K \delta (1 + \gamma) \leq 2 \kappa_0 (1+\gamma) c_2 \rho^{\frac{1}{\mu+1}} \delta^{\frac{\mu}{\mu+1}}.\]
In total we can conclude:
\[\|x^\dag- x_{n,K}^\delta\| \leq \left(\kappa_0^{\frac{1}{\mu+1}} (\tau + c_1)^{\frac{\mu}{\mu+1}} + 2 \kappa_0 (1+\gamma) c_2 \right) \rho^{\frac{1}{\mu+1}}\delta^{\frac{\mu}{\mu+1}}.\]
\qed

Regarding the computational expenses of Algorithm \ref{algorithm1}, we get the following result. 

\begin{theorem} \label{teorem3}
Let $A \in {\cal H}^r$, $x^\dag \in M_{\mu,\rho}(A)$, $0 < \mu \leq \mu_0 - 1$, $\mu_0 \geq 2$ and $x_{n,k}^\delta$ be computed according to Algorithm \ref{algorithm1}. Further, we assume
that the level of noise satisfies $\delta < \|f\|$ and that Algorithm \ref{algorithm1} is not stopped in the first iteration of the while loop. Then, the discretization level $n$ of the
solution $x_{n,K}^\delta$ is bounded by. 
\begin{equation} \label{inequalitydiscretizationlevel}
n < c_4 + c_5 \ln \frac{\rho}{\delta}
\end{equation}
with $c_4 = \frac{1}{r \ln 2} \ln \frac{c_2}{\gamma} \frac{2^{r+1} (1+2^{r+3})}{2^r-1}$ and $c_5 = \frac{1}{r \ln 2}  \frac{\mu +2 }{ \mu+1}$. Further, the number of 
calculated inner products $\langle Ae_i, e_j\rangle$ for the domain $\Gamma_n$ can be estimated as 
\begin{equation} \label{inequalityinformational}
\# \Gamma_n = 2^{2n}(n+1) <  c_3 \left( \frac{\rho}{\delta} \right)^{\frac{\mu+2}{r(\mu+1)}} \left( 1+ c_4 + c_5 \ln \frac{\rho}{\delta} \right)^{1 + \frac1r},
\end{equation}
with $c_3 = \left( \frac{c_2}{\gamma} 2^{2r+1} (1+2^{r+3}) \right)^{\frac1r}$,
\end{theorem}

{\bf Proof.} Since Algorithm \ref{algorithm1} is not stopped in the first iteration we have by definition of the index $K_{n-1}$ in \eqref{equation-n}:
        \begin{equation*} 
                (1+2^{r+3})2^{-2r (n-1)} (n-1) > \frac{\gamma \delta}{2 (K_{n-1}+1) \rho}.
        \end{equation*} 
This can be formulated equivalently as
        \begin{equation*} 
                \frac{2^{2r n}}{n-1} < 2^{2r+1} (1+2^{r+3}) \frac{\rho}{ \gamma \delta } (K_{n-1}+1).
        \end{equation*} 
By inequality \eqref{inequalitymaximalindex} in Lemma \ref{lemma5} we get now the estimate
        \begin{equation} \label{inequality2} 
                \frac{2^{2r n}}{n-1} < \frac{c_2}{ \gamma} 2^{2r+1} (1+2^{r+3}) \left( \frac{\rho}{\delta} \right)^{\frac{\mu+2}{\mu+1}}. 
        \end{equation} 
The Bernoulli inequality $1+(n-1)(2^r-1) \leq 2^{r(n-1)}$ yields the bound $ (n-1) \leq \frac{2^{rn-r}}{2^r-1}$. Thus, we obtain from \eqref{inequality2} the inequality
        \begin{equation*} 
                2^{r n} < \frac{c_2}{ \gamma} \frac{2^{r+1}}{2^r-1} (1+2^{r+3}) \left( \frac{\rho}{\delta} \right)^{\frac{\mu+2}{\mu+1}}. 
        \end{equation*} 
Now, taking the logarithm on both sides gives the desired inequality \eqref{inequalitydiscretizationlevel}. Finally, using again \eqref{inequality2} and \eqref{inequalitydiscretizationlevel} we obtain
\[ \# \Gamma_n = 2^{2n}(n+1) < \left(\frac{2^{2r n}}{n-1}\right)^{\frac1r} (n+1)^{1+\frac1r} < 
c_3\left( \frac{\rho}{\delta} \right)^{\frac{\mu+2}{r(\mu+1)}} \left( 1+ c_4 + c_5 \ln \frac{\rho}{\delta} \right)^{1 + \frac1r},\]
with $c_3 = \left( \frac{c_2}{\gamma} 2^{2r+1} (1+2^{r+3}) \right)^{\frac1r}$.
\qed

\begin{remark}
For a standard nonadaptive Galerkin scheme with quadratic domain $\Omega = [1, \ldots 2^{2n}]^2$ the computational expenses for a suitable large discretization level $n$ turn out to be of order
$\mathrm{O}((\frac{\rho}{\delta})^{\frac{2}{r}})$. This asymptotic result can be deduced from the error bounds given in \cite{PlVai} and the fact that $A \in H_r$. In comparison, by Theorem \ref{teorem3}
the computational costs of Algorithm \ref{algorithm1} are of order $\mathrm{O}\left((\frac{\rho}{\delta})^{\frac{\mu+2}{(\mu+1)r}} \left(\ln \frac{\rho}{\delta} \right)^{1+\frac1r} \right)$. In this sense, the adaptive algorithm presented 
in this article is more economic than regularization schemes using the standard Galerkin scheme. The same order of complexity was also 
shown for the adaptive scheme in \cite{MaaPerRamSol} using a Tikhonov regularization and in \cite{SolVol} for the regularization with the $1$-method (with different proofs). 
\end{remark}

\subsection{The balancing principle as stopping rule}

In this section, we consider a second adaptive algorithm to solve (\ref{eq_main}) combining again a method from the class $R_{n}^{\mu_0}$ $\mu_0 \geq 2$ with 
an adaptive discretization strategy. However, this time we use the balancing principle (see \eqref{Balanc_princp} in Algorithm \ref{algorithm2}) as stopping rule.

\begin{algorithm}
\caption{Adaptive algorithm to solve (\ref{eq_main}) using the balancing principle}
\label{algorithm2}

\begin{algorithmic}[H]
\vspace{2mm}
\STATE given data $A\in {\cal H}^r, \delta, f_\delta, \rho$.
\STATE choose control parameters $\gamma > 0$, $K_{sec} \in \Nn$.
\STATE choose discretization level $n$ such that $(1+2^{r+3})2^{-2rn} n < \frac{\gamma \delta}{2 \rho}$ holds. 
\STATE compute Galerkin information:
        $$
        (f_\delta, e_j), \quad j\in [1, 2^{2n}]; \qquad
        (Ae_i, e_j), \quad (i,j) \in \Gamma_{n}.
        $$
\WHILE {(balancing principle $=$ false)}  \vspace{2mm}

\STATE choose $K_n \in \Nn$ as maximal integer such that 
        \begin{equation} \label{equation-n2}
                (1+2^{r+3})2^{-2rn} n < \frac{\gamma \delta}{2 K_n \rho}
        \end{equation} 
        is satisfied. \\[2mm]
\FOR {($k = 1:K_n + K_{sec}$)} \vspace{2mm}
\STATE compute iterates of semiiterative method in the class $R_n^{\mu_0}$, $\mu_0 \geq 2$ (cf. \eqref{SimiItMeth}):
\STATE \[ x_{n,k}^\delta = x_{n,k-1}^\delta + ((1-\alpha_k) \omega_{k}-1) ( x_{n,k-1}^\delta - x_{n,k-2}^\delta) + 2 \omega_{k}  A_n^*( f_\delta - A_n x_{n,k-1}^\delta). \]
\ENDFOR\\[2mm]
\STATE compute the set
            \begin{equation}\label{Balanc_princp}
                D_n^+ = \{k: k \leq K_n, \|x_{n,k}^\delta - x_{n,j}^\delta\|\leq 8 (1 + \gamma) \kappa_0 j \delta \quad \text{for all $j$ with $k < j \leq K_n+K_{sec}$} \}.
            \end{equation}
\IF {$D_n^+=\varnothing$} 
\STATE increase discretization level $n \to n+1$.  
\ELSE 
\STATE balancing principle $=$ true 
\RETURN stopping index $K=\min\{k: k \in D_n^+\}$, discretization level $n$ and solution $x_{n,K}^\delta$ \\[2mm]
\ENDIF
\STATE compute new Galerkin information:
        $$
        (f_\delta, e_j), \quad j\in [2^{2n-2}, 2^{2n}]; \qquad
        (Ae_i, e_j), \quad (i,j) \in \Gamma_{n} \setminus \Gamma_{n-1}.
        $$
\ENDWHILE
\end{algorithmic}
\end{algorithm}

For technical purposes we need the index $K_{opt} :=\left\lceil\left(\frac{2 (1+\gamma) \delta}{\kappa_{\mu}\rho}\right)^{-\frac{1}{\mu+1}}\right\rceil$. It
is easy to see that $\kappa_{\mu} \rho K^{-\mu}_{opt} \leq 2(1+\gamma) K_{opt}\delta$ is satisfied and that for all $k\geq K_{opt}$ we have the inequality
\begin{equation}\label{k1}
 \kappa_{\mu} \rho k^{-\mu}\leq 2(1+\gamma) k \delta.
\end{equation}

\begin{theorem}\label{teorem2}
Let $A \in {\cal H}^r$ and $\mu_0 \geq 2$ for the class $R_n^{\mu_0}$. Then, Algorithm \ref{algorithm2} gives
an order optimal regularization method for the solution $x^\dag$ of \eqref{eq_main} in the class $M_{\mu,\rho}(A)$ for all $0 < \mu \leq \mu_0$. In particular,
the approximative solution $x_{n,K}^\delta$ given by Algorithm \ref{algorithm2} satisfies 
\begin{equation} \label{orderoptimalitybalancing}
\|x^\dag-x_{n,K}^\delta\|\leq C \rho^{\frac{1}{\mu+1}} \delta^{\frac{\mu}{\mu+1}},
\end{equation}
with $C = 12 \kappa_\mu^{\frac{1}{\mu+1}} \left(2(1+\gamma)\kappa_0 \right)^{\frac{\mu}{\mu+1}}$.
\end{theorem}

{\bf Proof.} 
We check that $K_{opt} \geq K$. For all $k\geq 1$, Corollary \ref{corollary-discretization} implies
$$
\| x^\delta_{n,k} - x^\delta_{n,K_{opt}}\| \leq \|x^\dag - x^\delta_{n,K_{opt}}\|+\|x^\dag - x^\delta_{n,k}\|
\leq \kappa_{\mu} \rho K^{-\mu}_{opt}+ 2(1+\gamma) \kappa_0 K_{opt} \delta + \kappa_{\mu} \rho k^{-\mu}+ 2(1+\gamma) \kappa_0 k \delta.
$$
Now, using (\ref{k1}) we get for all $k\geq K_{opt}$:
$$
\| x^\delta_{n,k} - x^\delta_{n,K_{opt}}\| \leq 4 (1+\gamma) \kappa_0 k \delta + 4(1+\gamma) \kappa_0 K_{opt}\delta \leq 8(1+\gamma) \kappa_0 k \delta.
$$
Thus, $K_{opt}\in D_n^+$ if $K_{opt} \leq K_n$ and $K_{opt} \geq K.$
Therefore, by the balancing principle \eqref{Balanc_princp} the total error can be bounded as follows:
\begin{align*}
 \|x^\dag- x^\delta_{n,K}\| &\leq\|x^\dag- x^\delta_{n,K_{opt}}\|+\| x^\delta_{n,K_{opt}} - x^\delta_{n,K}\|\leq 
 \kappa_{\mu} \rho K^{-\mu}_{opt} + 2(1+\gamma) \kappa_0 K_{opt} \delta + 8(1+\gamma) \kappa_0 K_{opt}\delta \\
 & \leq 12 (1 + \gamma) \kappa_0 K_{opt}\delta = 12 (1 + \gamma) \kappa_0 \left\lceil \left(\frac{2(1+\gamma) \kappa_0 \delta}{\kappa_{\mu}\rho}\right)^{-\frac{1}{\mu+1}}\right\rceil \delta  \\
 & \leq 24 (1+\gamma) \kappa_0 \left(\frac{2(1+\gamma) \kappa_0 \delta}{\kappa_{\mu}\rho}\right)^{-\frac{1}{\mu+1}}\delta= C \rho^\frac{1}{\mu+1} \delta^\frac{\mu}{\mu+1}.
\end{align*}
Thus, Algorithm \ref{algorithm2} yields an order optimal reconstruction scheme for all $0 < \mu \leq \mu_0$. It follows again from \cite[Theorem 2.1]{Plato} that
Algorithm \ref{algorithm2} is also a regularization method for \eqref{eq_main}. \qed

\section{Examples and numerical tests}

In this final section, we present a simple test equation, in which the preliminary assumptions of Theorem \ref{teorem1} and \ref{teorem2} are satisfied. 
With help of this example we test the convergence order and the performance of the introduced adaptive algorithms.  
As a simple example of a linear problem (\ref{eq_main}) in which the operator $A$ is in the class ${\cal H}^2$, 
we consider in $X = L^2([0,1])$ the Fredholm integral equation of the first kind (see \cite[Example 12.4.1.]{DelvesMohamed}
\begin{equation}
A x(t) = \int_0^1 k(s,t) x(s) ds = f(t),\qquad 0 \leq t \leq 1,  \label{equation-generaltestequation}
\end{equation}
with the kernel
\[ 
k(s,t) = \left\{ \begin{array}{ll} t(s-1) & 0\leq t < s \leq 1, \\ s(t-1) & 0\leq s \leq t \leq 1. \end{array} \right. 
\]
The self-adjoint operator $A$ maps $L^2([0,1])$ into the Sobolev space $W^{2,2}([0,1])$. In the setting of boundary value problems, the application
of $A$ to $x$ corresponds to the solution of the boundary value problem $f''(t) = x(t)$ with homogeneous boundary conditions $f(0) = f(1) = 0$. It is well-known that
the functions 
\[e_k(t) = \sqrt{2}\sin (\pi k t), \quad k \in \Nn, \; t \in [0,1],\]
form an orthonormal basis of eigenfunctions of the operator $A$ with corresponding eigenvalues $\lambda_k = - (\pi k)^{-2}$. Therefore, we have 
$\|A\| \leq \pi^{-2}$ and $\|( I - P_m) A \| = \|A ( I - P_m)\| \leq (\pi (m+1))^{-2}$. This implies in particular $\pi^2 A \in {\cal H}^2$. 

In \eqref{equation-generaltestequation}, we consider the two different right hand sides $f_1$ and $f_2$ given by
\begin{align*}
f_1(t) &= t^3(1-t)^3, \\
f_2(t) &= \textstyle \frac{t^3}3 - \max\left(0,t-\frac12 \right)^2 - \frac{t}{12}.
\end{align*}
The corresponding exact solutions of the inverse problem \eqref{equation-generaltestequation} are given by
\begin{align*}
x_1^\dag(t) &= 18 t(5t^3-10t^2+6t-1) \quad \text{with $x_1^\dag \in M_{\mu,\rho}$ for $0 < \mu < 1.25$}, \\
x_2^\dag(t) &= 2t - \mathrm{sign}(2t-1) - 1 \quad \text{with $x_2^\dag \in M_{\mu,\rho}$ for $0 < \mu < 0.25$}.
\end{align*}

\subsection{Test of Algorithm 1}
\begin{figure}[H]  \caption{The error and the stopping index $K$ of Algorithm \ref{algorithm1} compared 
to the expected error $\delta^{\frac{1.25}{1+1.25}}$ and the expected stopping index $\delta^{-\frac{1}{1+1.25}}$
of the first test problem $A x_1 = f_1$.}\label{figure-1} 
  \begin{minipage}{0.5\textwidth}
  \centering
  \includegraphics[width=\textwidth]{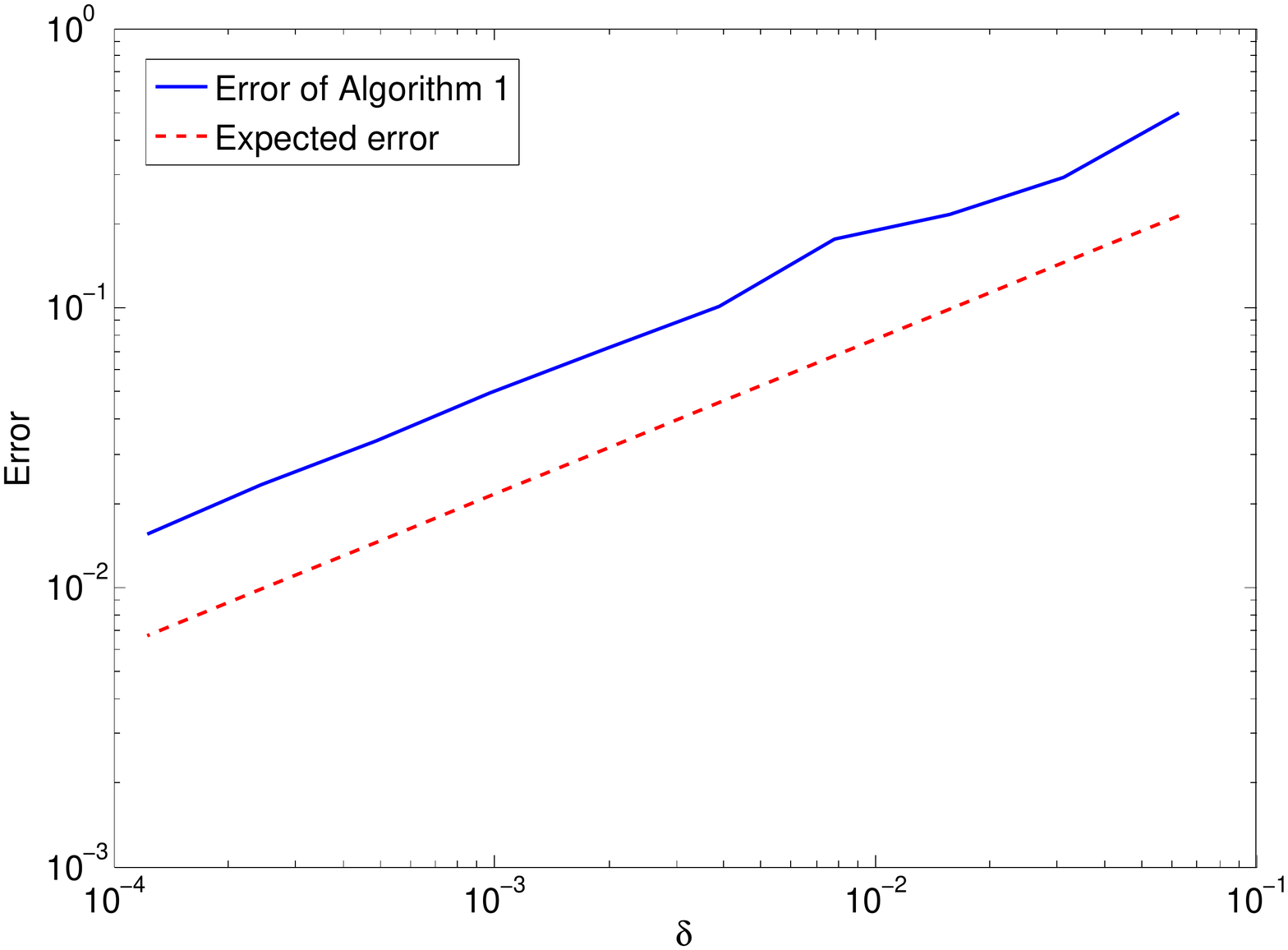}
  \end{minipage}\hfill
  \begin{minipage}{0.5\textwidth}
  \centering
  \includegraphics[width=\textwidth]{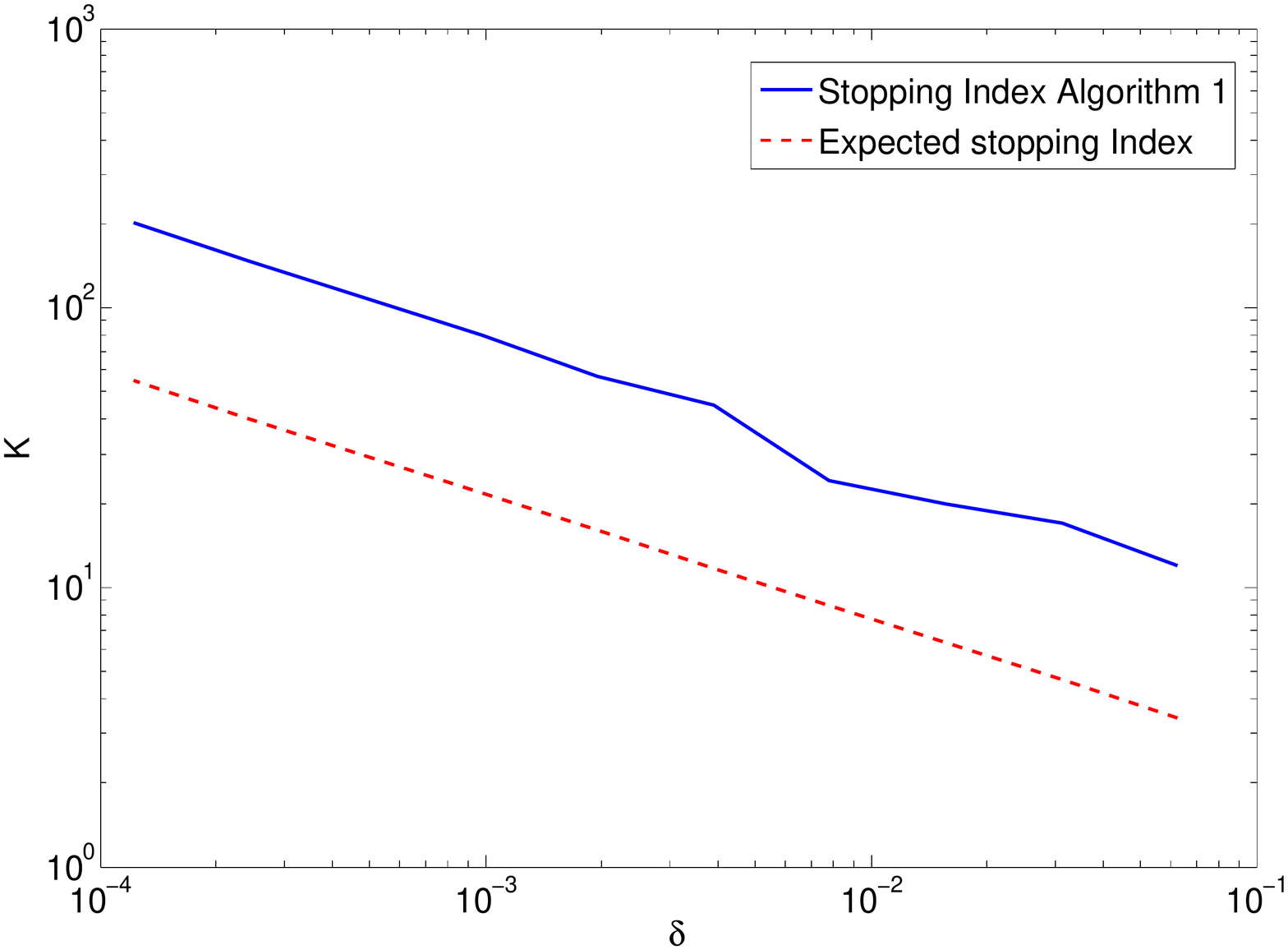}
  \end{minipage}
\end{figure}

\begin{figure}[H] \caption{The error and the stopping index $K$ of Algorithm \ref{algorithm1} compared 
to the expected error $\delta^{\frac{0.25}{1+0.25}}$ and the expected stopping index $\delta^{-\frac{1}{1+0.25}}$
of the second test problem $A x_2 = f_2$.} \label{figure-2} 
  \begin{minipage}{0.5\textwidth}
  \centering
  \includegraphics[width=\textwidth]{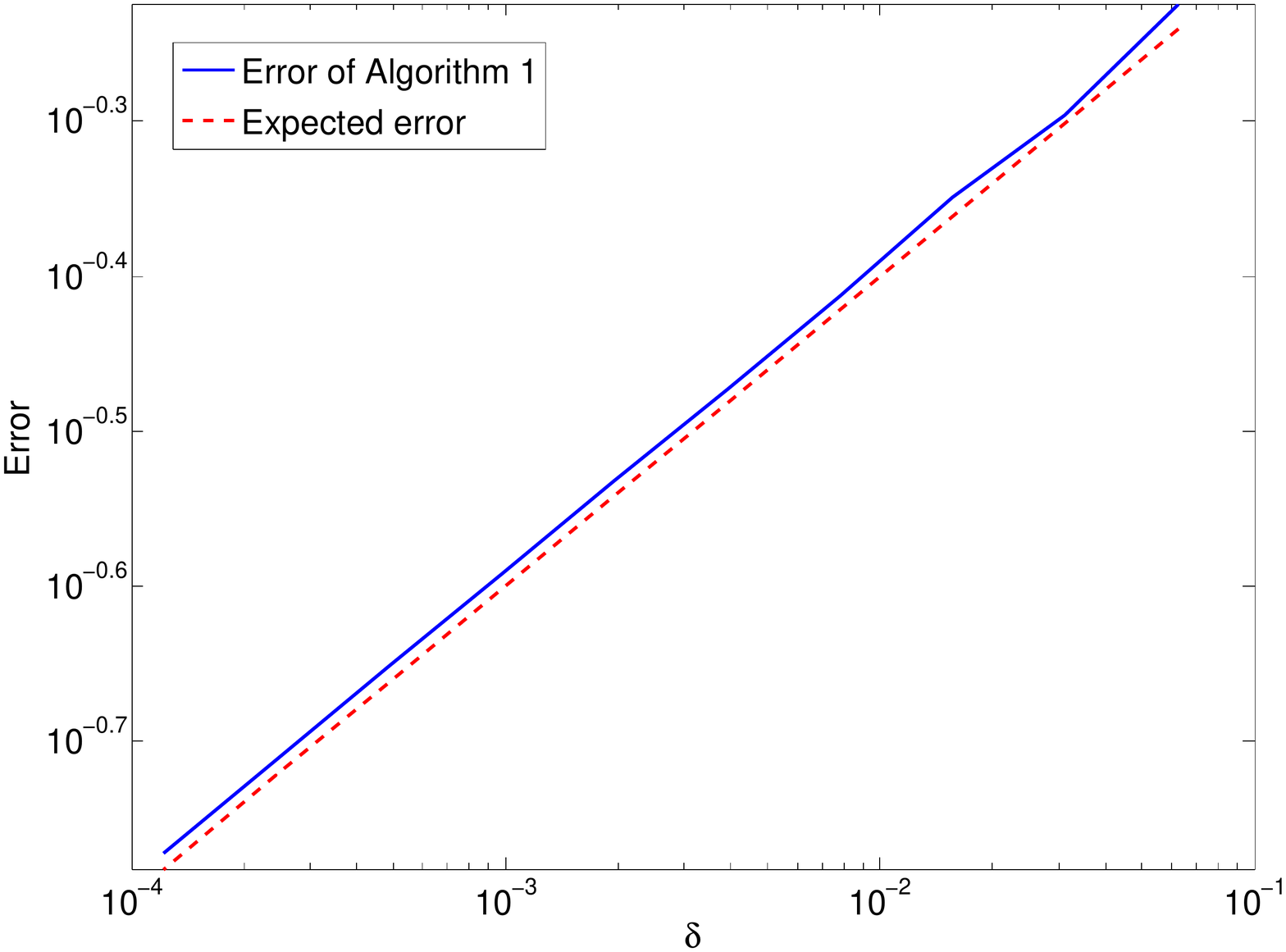} 
  \end{minipage}\hfill
  \begin{minipage}{0.5\textwidth}
  \centering
  \includegraphics[width=\textwidth]{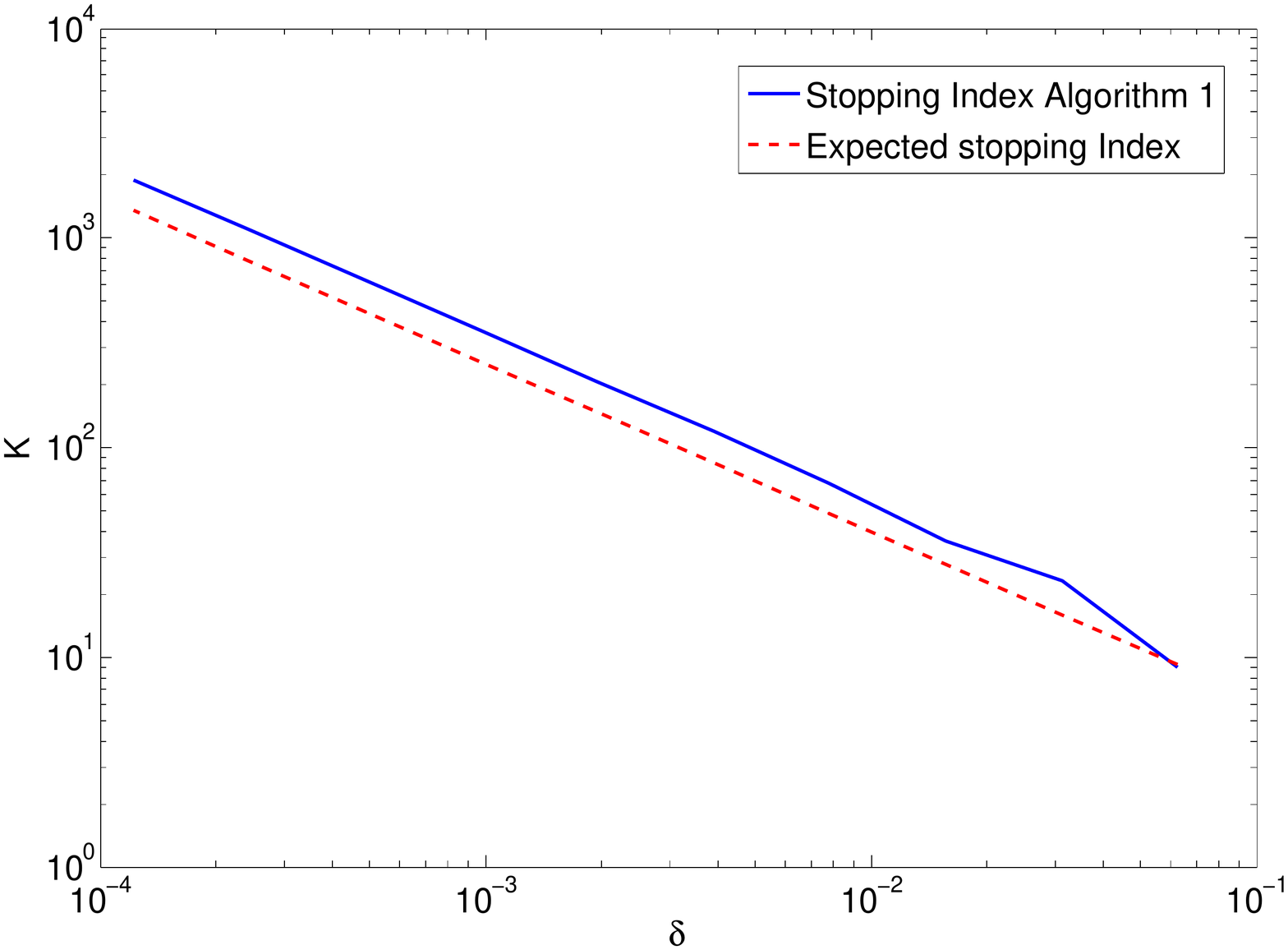}
  \end{minipage}
\end{figure}

We test first Algorithm 1 for the two right hand sides $f_1$ and $f_2$. As semiiterative method we chose the $\nu$-method with $\nu = 1.5$. Since the qualification of this method is $\mu_0 = 3$ it can
be applied in Algorithm 1 for both test examples. In this way we have $\kappa_0 = 1$ and $\kappa_2 = 6$. Choosing $\gamma = \frac12$, we can take $\tau = 1.01 + \sqrt{\frac{13}{8}}$ as parameter for
the discrepancy principle. We set $\rho = 1$ and generate perturbed right hand sides $f_{i,\delta}, i = 1,2$ for different values of $\delta > 0$. Now, we use Algorithm 1 to compute the approximate solutions 
$x_{i,n,K}^\delta$ of $x_i^\dag$. The resulting errors and stopping indices of Algorithm \ref{algorithm1} are displayed in Table \ref{table-1}, \ref{table-2} and Figure \ref{figure-1}, \ref{figure-2}

\begin{table}[H] \caption{Results of Algorithm \ref{algorithm1} to solve the test problem $A x_1 = f_1$.
}\label{table-1} \centering
  \begin{tabular}[t]{|ccccccc|} \hline 
  $\nu$ & $\delta$ & $n$ & $K_n$ & $K$ & $\frac{\|x_{1,n,K}^\delta -x_1^\dag\|}{\|x_1^\dag\|}$ & $\delta^{\frac{1.25}{1+1.25}}$ \\[2mm] \hline 
 1.5  &  0.062500  &   6  &     37  &    12  &  0.49975111  &  0.21431100 \\
 1.5  &  0.031250  &   6  &     19  &    17  &  0.29238913  &  0.14581613 \\
 1.5  &  0.015625  &   7  &    125  &    20  &  0.21650878  &  0.09921257 \\
 1.5  &  0.007812  &   7  &     63  &    24  &  0.17715080  &  0.06750373 \\
 1.5  &  0.003906  &   8  &    435  &    45  &  0.10086226  &  0.04592920 \\
 1.5  &  0.001953  &   8  &    218  &    57  &  0.07100275  &  0.03125000 \\
 1.5  &  0.000977  &   8  &    109  &    80  &  0.04971398  &  0.02126234 \\
 1.5  &  0.000488  &   9  &    774  &   108  &  0.03362040  &  0.01446679 \\
 1.5  &  0.000244  &   9  &    387  &   147  &  0.02322422  &  0.00984313 \\
 1.5  &  0.000122  &  10  &   2784  &   203  &  0.01549616  &  0.00669722 \\ \hline
  \end{tabular}
\end{table}

% The errors of Algorithm \ref{algorithm1} for the right hand side $f_2 $ are displayed in Table \ref{table-2} and Figure \ref{figure-2}

\begin{table}[H] \caption{Results of Algorithm \ref{algorithm1} to solve the test problem $A x_2 = f_2$.
}\label{table-2} \centering
  \begin{tabular}[t]{|ccccccc|} \hline 
  $\nu$ & $\delta$ & $n$ & $K_n$ & $K$ & $\frac{\|x_{2,n,K}^\delta -x_2^\dag\|}{\|x_2^\dag\|}$ & $\delta^{\frac{0.25}{1+0.25}}$ \\[2mm] \hline 
 1.5  &  0.062500  &   6  &     16  &     9  &  0.59696031  &  0.57434918 \\
 1.5  &  0.031250  &   7  &    105  &    23  &  0.50523819  &  0.50000000 \\
 1.5  &  0.015625  &   7  &     53  &    36  &  0.44800149  &  0.43527528 \\
 1.5  &  0.007812  &   8  &    366  &    68  &  0.38638037  &  0.37892914 \\
 1.5  &  0.003906  &   8  &    183  &   120  &  0.33629235  &  0.32987698 \\
 1.5  &  0.001953  &   9  &   1300  &   207  &  0.29364826  &  0.28717459 \\
 1.5  &  0.000977  &   9  &    650  &   361  &  0.25566471  &  0.25000000 \\
 1.5  &  0.000488  &  10  &   4678  &   625  &  0.22295988  &  0.21763764 \\
 1.5  &  0.000244  &  10  &   2339  &  1091  &  0.19402742  &  0.18946457 \\
 1.5  &  0.000122  &  11  &  17010  &  1901  &  0.16890368  &  0.16493849 \\ \hline
  \end{tabular}
\end{table}

\subsection{Test of Algorithm 2}
\begin{figure}[H] \caption{The error and the stopping index $K$ of Algorithm \ref{algorithm2} compared 
to the expected error $\delta^{\frac{1.25}{1+1.25}}$ and the expected stopping index $\delta^{-\frac{1}{1+1.25}}$
of the first test problem $A x_1 = f_1$.} \label{figure-3} 
  \begin{minipage}{0.5\textwidth}
  \centering
  \includegraphics[width=\textwidth]{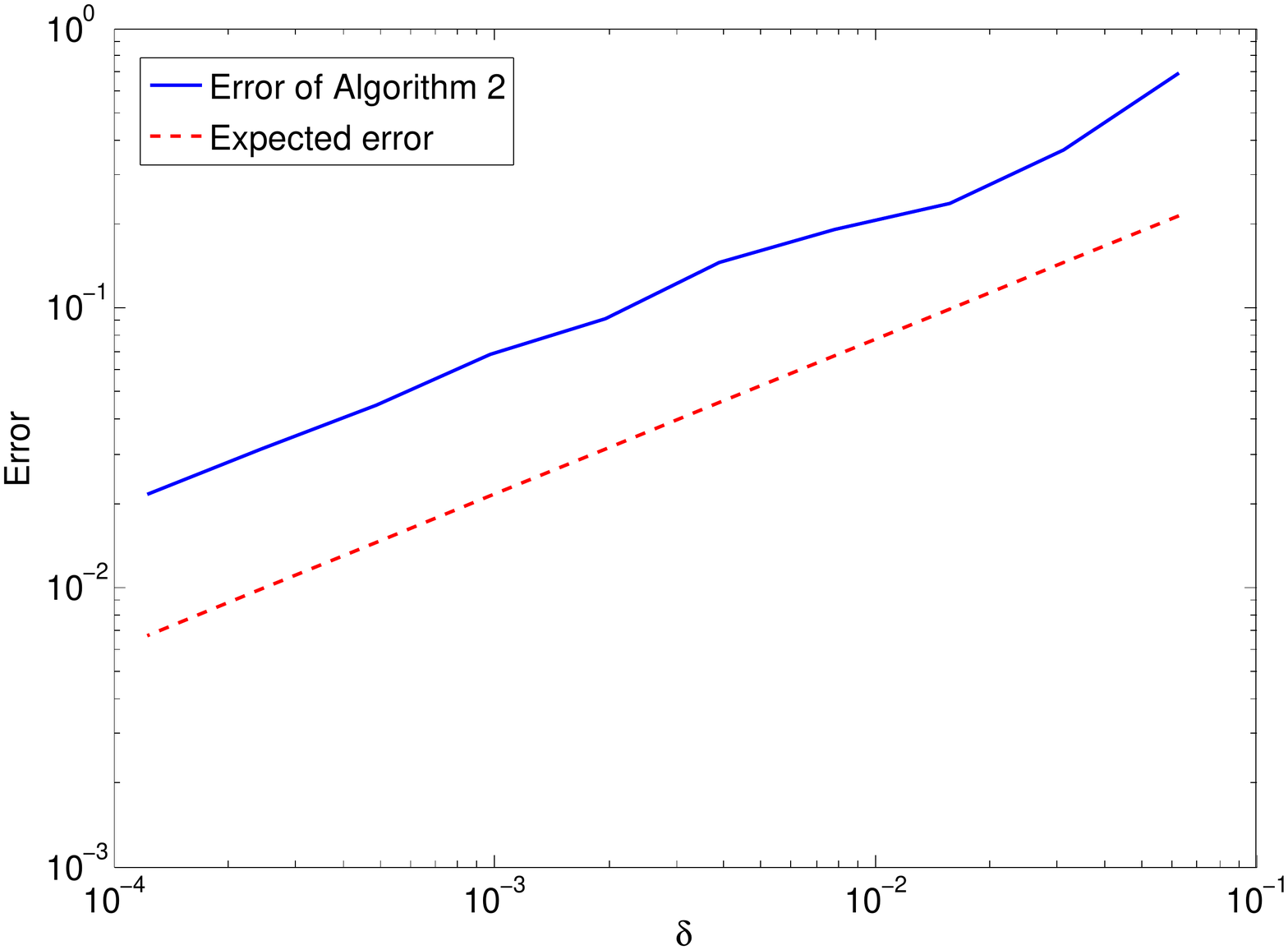}
  \end{minipage}\hfill
  \begin{minipage}{0.5\textwidth}
  \centering
  \includegraphics[width=\textwidth]{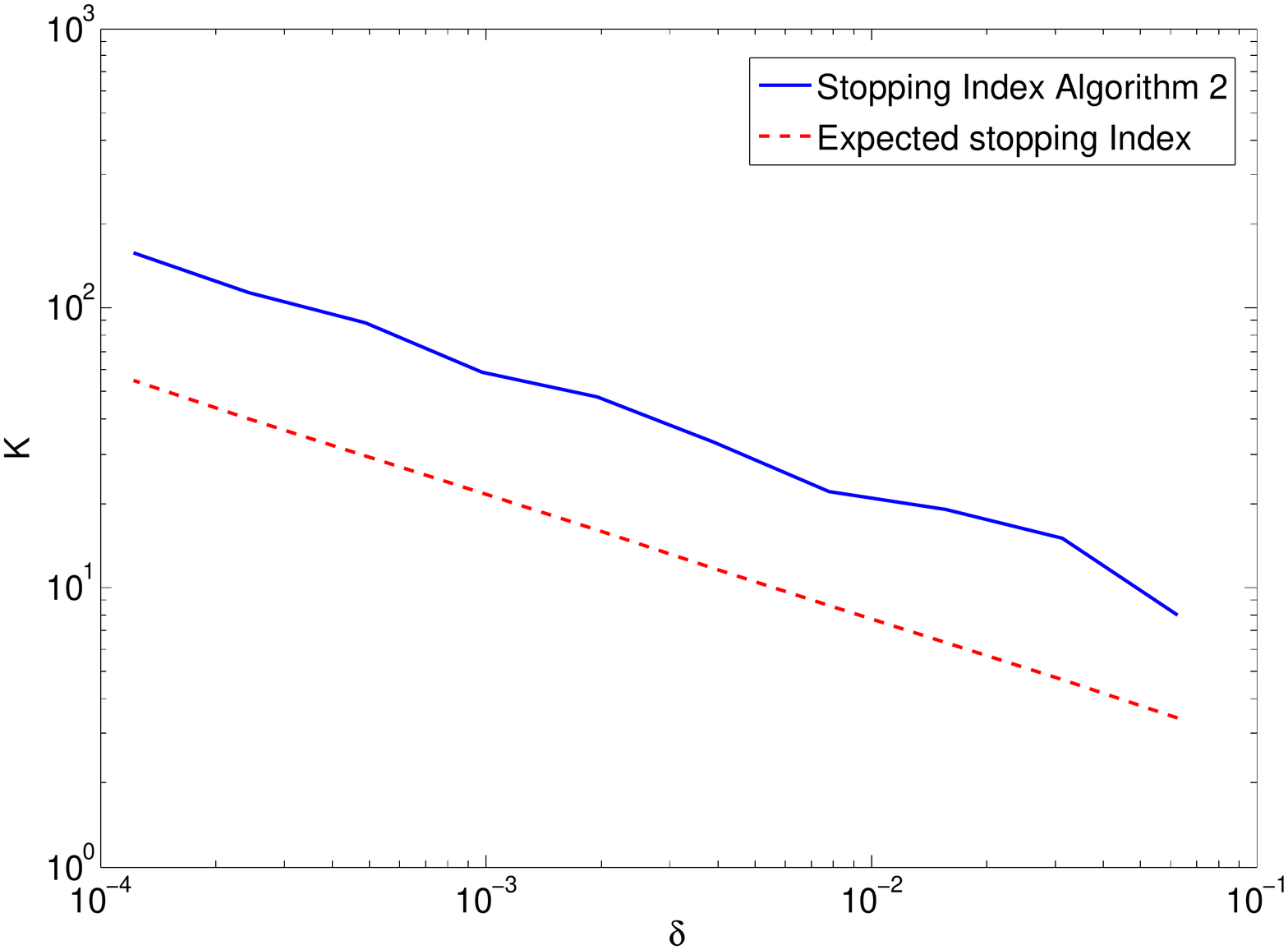}
  \end{minipage}
\end{figure}

\begin{figure}[H]  \caption{The error and the stopping index $K$ of Algorithm \ref{algorithm2} compared 
to the expected error $\delta^{\frac{0.25}{1+0.25}}$ and the expected stopping index $\delta^{-\frac{1}{1+0.25}}$
of the first test problem $A x_2 = f_2$.} \label{figure-4} 
  \begin{minipage}{0.5\textwidth}
  \centering
  \includegraphics[width=\textwidth]{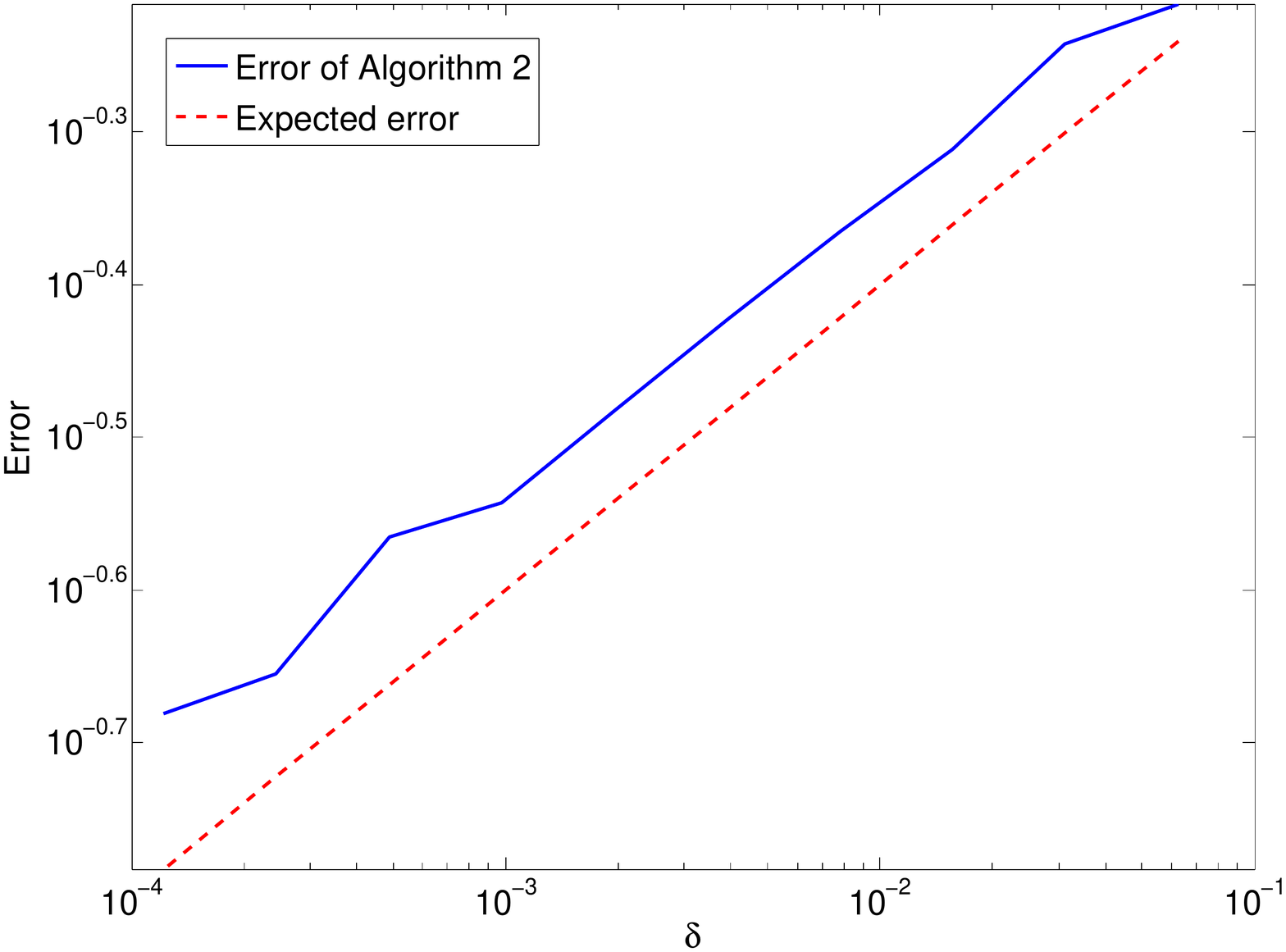}
  \end{minipage}\hfill
  \begin{minipage}{0.5\textwidth}
  \centering
  \includegraphics[width=\textwidth]{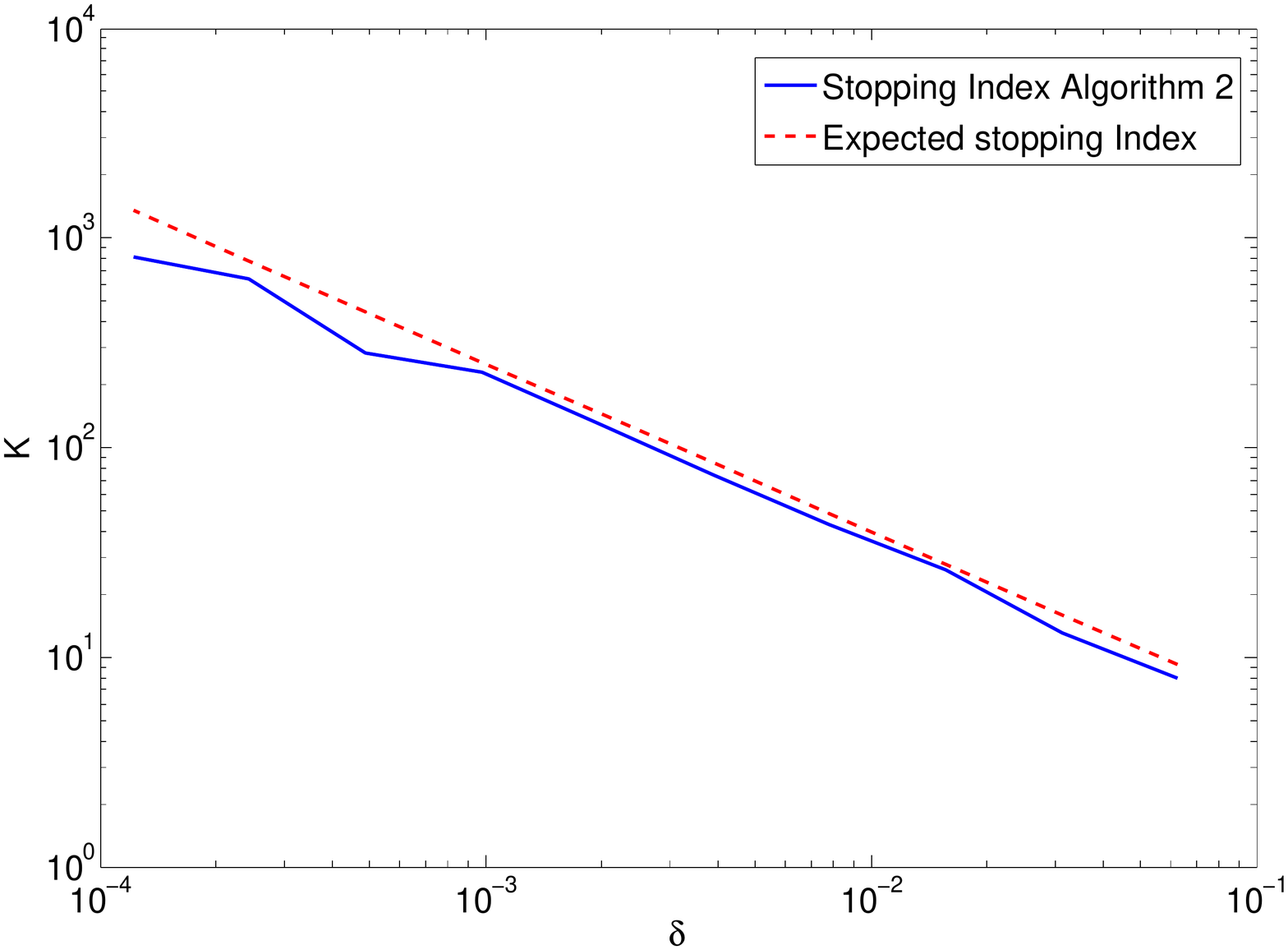}
  \end{minipage}
\end{figure}

Now we test Algorithm \ref{algorithm2} for the two right hand sides $f_1$ and $f_2$. As semiiterative method we choose again the $\nu$-method with $\nu = 1.5$. 
As further control parameter, we choose $\gamma = \frac12$. Also for Algorithm \ref{algorithm2}, we set $\rho = 1$ and generate perturbed right hand sides 
$f_{i,\delta}, i = 1,2$, for different values of $\delta > 0$. With Algorithm \ref{algorithm2} we then compute approximate solutions $x_{i,n,K}^\delta$ of $x_i^\dag$. 
The resulting errors for the two test problems are displayed in Table \ref{table-3}, \ref{table-4} and Figure \ref{figure-3}, \ref{figure-4}.

\begin{table}[H] \caption{Results of Algorithm \ref{algorithm2} to solve the test problem $A x_1 = f_1$.
}\label{table-3} \centering
\begin{tabular}[t]{|ccccccc|} \hline 
  $\nu$ & $\delta$ & $n$ & $K_n$ & $K$ & $\frac{\|x_{1,n,K}^\delta -x_1^\dag\|}{\|x_1^\dag\|}$ & $\delta^{\frac{1.25}{1+1.25}}$ \\[2mm] \hline 
 1.5  &  0.062500  &   6  &     37  &     8  &  0.68979661  &  0.21431100 \\
 1.5  &  0.031250  &   6  &     19  &    15  &  0.36601474  &  0.14581613 \\
 1.5  &  0.015625  &   7  &    125  &    19  &  0.23679445  &  0.09921257 \\
 1.5  &  0.007812  &   7  &     63  &    22  &  0.18993287  &  0.06750373 \\
 1.5  &  0.003906  &   8  &    435  &    33  &  0.14615533  &  0.04592920 \\
 1.5  &  0.001953  &   8  &    218  &    48  &  0.09181469  &  0.03125000 \\
 1.5  &  0.000977  &   8  &    109  &    59  &  0.06784866  &  0.02126234 \\
 1.5  &  0.000488  &   9  &    774  &    88  &  0.04481807  &  0.01446679 \\
 1.5  &  0.000244  &   9  &    387  &   114  &  0.03125762  &  0.00984313 \\
 1.5  &  0.000122  &   9  &    194  &   158  &  0.02144644  &  0.00669722 \\ \hline
\end{tabular}
\end{table}

\begin{table}[H] \caption{Results of Algorithm \ref{algorithm2} to solve the test problem $A x_2 = f_2$.
}\label{table-4} \centering
\begin{tabular}[t]{|ccccccc|} \hline 
  $\nu$ & $\delta$ & $n$ & $K_n$ & $K$ & $\frac{\|x_{2,n,K}^\delta -x_2^\dag\|}{\|x_2^\dag\|}$ & $\delta^{\frac{0.25}{1+0.25}}$ \\[2mm] \hline 
 1.5  &  0.062500  &   6  &     16  &     8  &  0.60790728  &  0.57434918 \\
 1.5  &  0.031250  &   7  &    105  &    13  &  0.57256321  &  0.50000000 \\
 1.5  &  0.015625  &   7  &     53  &    26  &  0.48809868  &  0.43527528 \\
 1.5  &  0.007812  &   8  &    366  &    43  &  0.43126730  &  0.37892914 \\
 1.5  &  0.003906  &   8  &    183  &    74  &  0.37805221  &  0.32987698 \\
 1.5  &  0.001953  &   9  &   1300  &   131  &  0.32892273  &  0.28717459 \\
 1.5  &  0.000977  &   9  &    650  &   228  &  0.28663978  &  0.25000000 \\
 1.5  &  0.000488  &   9  &    325  &   281  &  0.27204892  &  0.21763764 \\
 1.5  &  0.000244  &  10  &   2339  &   643  &  0.22139020  &  0.18946457 \\
 1.5  &  0.000122  &  10  &   1170  &   818  &  0.20848631  &  0.16493849 \\ \hline
\end{tabular}
\end{table}

The graphs in Table \ref{table-1} - \ref{table-4} confirm the 
theoretical results on the order optimality of Algorithm \ref{algorithm1} and \ref{algorithm2}. The parameters $\tau$ and
$\gamma$ in the two algorithms are chosen conservatively such that the assumptions in Theorem \ref{teorem1} and \ref{teorem2} are satisfied. Neglecting these theoretical preconditions on the parameters it is possible to further improve the numerical results. In particular for Algorithm \ref{algorithm1}, a smaller choice of $\tau$ yields better results for the approximation error $\frac{\|x_{n,K}^\delta -x^\dag\|}{\|x^\dag\|}$. Also, increasing the parameter $\gamma$ leads to smaller discretization levels $n$, as compared to the ones displayed in Table \ref{table-1} - \ref{table-4}, and reduces the numerical costs of the algorithms.

\end{document}